\documentclass{article}
\usepackage{amssymb}
\usepackage{graphicx}
\usepackage{amsmath}

\input{tcilatex}
\begin{document}

\begin{center}
\ \ \ \ \ \ \ \ \ \ \ \ \ \ \ \ \ \ \ \ \ \ \ \ \ \ \ \ \ \ \ \ \ \ \ \ \ \
\ \ \ \ \ \ \ \ 

\bigskip 

\bigskip \textbf{\ Differential operators in exterior domain and application}

\ \ \ \ \ \ \ \ \ \ \ \ \ \ \ \ \ \ \ \ \ 

\bigskip\ \ \ \ \ \ \ \ \ \ \ \ \ \ \ \ \ \ \textbf{Veli B. Shakhmurov}\ \ \ 

Okan University, Department of Mechanical Engineering, Akfirat, Tuzla 34959
Istanbul, Turkey, E-mail: veli.sahmurov@okan.edu.tr

\ \ \ \ \ \ \ \ \ \ \ \ \ \ \ \ \ \ \ \ \ \ \ 

\textbf{AMS: \ 35xx,\ 47Fxx, 47Hxx, 35Pxx}\ \ \ \ \ \ \ \ \ \ \ \ \ 

\textbf{Abstract}
\end{center}

The abstract elliptic and parabolic equations on exterior domain are
considered. The equations have top-order variable coefficients. The
separability properties of boundary value problems for elliptic equation and
well-posedness of the Cauchy problem for parabolic equations are
established. In application, the well-posedness of Wentzell-Robin type mixed
probem for parabolic equation, Cauchy problem for anisotropic parabolic
equations and system of parabolic equations are derived


\begin{center}
\ \ \textbf{1. Introduction, notations and background }
\end{center}

Boundary value problems (BVPs) for differential-operator equations (DOEs)
have been studied extensively by many researchers (see $\left[ \text{3, 5,
8-23, 26}\right] $ and the references therein). A comprehensive introduction
to the DOEs and historical references may be found in $\left[ 13\right] $
and $\left[ 26\right] .$ The maximal regularity properties for differential
operator equations have been studied in $\left[ \text{2}\right] ,$ $\left[ 
\text{8}\right] ,$ $\left[ 9\right] $ and $\left[ \text{17-23}\right] $ for
instance. The main objective of the present paper is to discuss the exterior
BVPs for the following DOE with variable coefficients 
\begin{equation}
\varepsilon au^{\left( 2\right) }\left( x\right) +Au\left( x\right)
+\varepsilon ^{\frac{1}{2}}A_{1}u^{\left( 1\right) }\left( x\right)
+A_{0}u\left( x\right) =f\left( x\right) ,\text{ }x\in \sigma ,  \tag{1.1}
\end{equation}%
\begin{equation*}
\sum\limits_{i=0}^{\mu _{1}}\alpha _{i}\varepsilon ^{\nu _{i}}u^{\left(
i\right) }\left( 0\right) =0,\text{ }\sum\limits_{i=0}^{\mu _{2}}\beta
_{i}\varepsilon ^{\nu _{i}}u^{\left( i\right) }\left( b\right) =0\text{, }
\end{equation*}%
where $\sigma $ is an exterier domain, i.e. $\sigma $ $=\left( -\infty
,\infty \right) /\left[ 0,b\right] ,$ $a=a\left( x\right) $ is a
complex-valued function, $\varepsilon $ is a positive parameter, $\nu _{i}=%
\frac{i}{2}+\frac{1}{2p},$ $p\in \left( 1,\infty \right) ;$ $A=A\left(
x\right) $, \ $A_{j}=A_{j}\left( x\right) $ are linear operator functions in
a Banach space $E$, $\alpha _{i}$, $\beta _{i}$ are complex numbers, $\mu
_{k}\in \left\{ 0,1\right\} .$

In this paper, the $E$-valued $L_{p}$-separability properties of this
problem is obtained. Especially, we prove that the corresponding
differential operator is $R$-positive and also is a negative generator of
the analytic semigroup.

Note that, the principal part of the corresponding differential operator is
non selfadjoint. Nevertheless, the sharp uniform coercive estimates for\ the
resolvent of corresponding differential operators are established. In
section 6, nonlocal BVP for degenerate abstract elliptic equation considered
in the moving domain. By using the maximal regularity properties of linear
problem $\left( 1.1\right) $ we derive the existence and uniqueness of BVP
for the following nonlinear degenerate abstract equation%
\begin{equation}
a\left( x\right) u^{\left( 2\right) }\left( x\right) +B\left( x,u,u^{\left(
1\right) }\right) u\left( x\right) =F\left( x,u,u^{\left( 1\right) }\right)
+f\left( x\right) ,  \tag{1.2}
\end{equation}%
in exterior domain, where $a$ is a complex valued function, $B$ and $F$ are
nonlinear operator in a Banach space $E.$

Then, by using the separability properties of the elliptic problem $\left(
1.1\right) $, the $L_{\mathbf{p}}\left( \sigma _{T};E\right) $
well-posedness is established for the following parabolic interior mixed
problem 
\begin{equation*}
\frac{\partial u}{\partial t}+\varepsilon a\frac{\partial ^{2}u}{\partial
x^{2}}+Au+\varepsilon ^{\frac{1}{2}}A_{1}\frac{\partial u}{\partial x}%
+A_{0}u=f\left( t,x\right) \text{, }t\in \left( 0,T\right) \text{, }x\in
\sigma ,
\end{equation*}%
\begin{equation}
\sum\limits_{i=0}^{\mu _{1}}\alpha _{i}\varepsilon ^{\nu _{i}}u^{\left(
i\right) }\left( t,0\right) =0,\text{ }\sum\limits_{i=0}^{\mu _{2}}\beta
_{i}\varepsilon ^{\nu _{i}}u^{\left( i\right) }\left( t,b\right) =0\text{, }
\tag{1.3}
\end{equation}%
\begin{equation*}
u\left( 0,x\right) =0\text{, }x\in \sigma .
\end{equation*}%
Here 
\begin{equation*}
\sigma _{T}=\sigma \times \left( 0,T\right) ,\text{ }\mathbf{p=}\left(
p_{1},p\right)
\end{equation*}%
and $L_{\mathbf{p}}\left( \sigma _{T};E\right) $ denotes the space of all $E$%
-valued $\mathbf{p}$-summable\ functions with mixed norm i.e., the space of
all $E$-valued measurable functions $f$ defined on $\sigma _{T}$ for which 
\begin{equation*}
\left\Vert f\right\Vert _{L_{\mathbf{p}}\left( \sigma _{T};E\right) }=\left(
\int\limits_{0}^{T}\left( \int\limits_{\sigma }\left\Vert f\left( t,x\right)
\right\Vert _{E}^{p_{1}}dt\right) ^{\frac{p}{p_{1}}}dx\right) ^{\frac{1}{p}%
}<\infty .
\end{equation*}

Moreover, let we choose $E=L_{2}\left( 0,1\right) $ in $\left( 1.1\right) $
and $A$ to be differential operator with generalized Wentzell-Robin boundary
condition defined by 
\begin{equation*}
D\left( A\right) =\left\{ u\in W_{p_{1}}^{2}\left( 0,1\right) ,\text{ }%
B_{j}u=Au\left( j\right) +\dsum\limits_{i=0}^{1}\alpha _{ij}u^{\left(
i\right) }\left( j\right) ,\text{ }j=0,1\right\} ,\text{ }
\end{equation*}%
\begin{equation*}
\text{ }Au=a_{1}u^{\left( 2\right) }+b_{1}u^{\left( 1\right) }+cu,
\end{equation*}%
where $\alpha _{ij}$ are complex numbers, $a_{1},$ $b_{1}$, $c$ are
complex-valued functions and $u^{\left( 0\right) }\left( x\right) =u\left(
x\right) $. Then, we get the $L_{\mathbf{\tilde{p}}}\left( \Omega \right) -$
well- posedness of the following Wentzell-Robin type mixed problem for
parabolic equation

\begin{equation}
\frac{\partial u}{\partial t}+\varepsilon a\frac{\partial ^{2}u}{\partial
x^{2}}+a_{1}\frac{\partial ^{2}u}{\partial y^{2}}+b_{1}\frac{\partial u}{%
\partial y}+cu=f\left( t,x,y\right) \text{, }  \tag{1.4}
\end{equation}%
\begin{equation*}
\sum\limits_{i=0}^{\mu _{1}}\alpha _{i}\varepsilon ^{\nu _{i}}u^{\left(
i\right) }\left( t,0\right) =0,\text{ }\sum\limits_{i=0}^{\mu _{2}}\beta
_{i}\varepsilon ^{\nu _{i}}u^{\left( i\right) }\left( t,b\right) =0\text{, }
\end{equation*}%
\ \ \ 
\begin{equation}
B_{j}u=0\text{, }j=0,1,\text{ }t\in \left( 0,T\right) \text{, }x\in \sigma ,%
\text{ }y\in \left( 0,1\right) ,  \tag{1.5}
\end{equation}

\begin{equation*}
u\left( 0,x,y\right) =0\text{, }x\in \sigma ,\text{ }y\in \left( 0,1\right) ,%
\text{ }
\end{equation*}%
where $\mathbf{\tilde{p}=}\left( \mathbf{p,}2\right) $, $\varepsilon $ is a
small parameter and $\Omega =\sigma _{T}\times \left( 0,1\right) .$

Note that, the regularity properties of Wentzell-Robin type BVP for elliptic
equations were studied e.g. in $\left[ \text{41, 42}\right] $ and the
references therein. The maximal regularity properties of DOEs in Banach
spaces were considered e.g. in $\left[ \text{2, 4, 9, 16, 21-23, 25}\right] $%
.

Let $L_{p}\left( \Omega ;E\right) $ denote the space of strongly measurable $%
E$-valued functions that are defined on $\Omega $ with the norm

\begin{equation*}
\left\Vert f\right\Vert _{L_{p}}=\left\Vert f\right\Vert _{L_{p}\left(
\Omega ;E\right) }=\left( \int \left\Vert f\left( x\right) \right\Vert
_{E}^{p}dx\right) ^{\frac{1}{p}},\text{ }1\leq p<\infty .
\end{equation*}

\ The Banach space\ $E$ is called an $UMD$-space if\ the Hilbert operator $%
\left( Hf\right) \left( x\right) =\lim\limits_{\varepsilon \rightarrow
0}\int\limits_{\left\vert x-y\right\vert >\varepsilon }\frac{f\left(
y\right) }{x-y}dy$ \ is bounded in $L_{p}\left( \mathbb{R},E\right) ,$ $p\in
\left( 1,\infty \right) $ (see. e.g. $\left[ 7\right] $). $UMD$ spaces
include e.g. $L_{p}$, $l_{p}$ spaces and Lorentz spaces $L_{pq},$ $p$, $q\in
\left( 1,\infty \right) $.

Let $\mathbb{R}$ denote the set of real numbers, $\mathbb{C}$ be the set of
the complex numbers and\ 
\begin{equation*}
S_{\varphi }=\left\{ \lambda ;\text{ \ }\lambda \in \mathbb{C}\text{, }%
\left\vert \arg \lambda \right\vert \leq \varphi \right\} \cup \left\{
0\right\} ,\text{ }0\leq \varphi <\pi .
\end{equation*}

\ Let $E_{1}$ and $E_{2}$ be two Banach spaces. $L\left( E_{1},E_{2}\right) $
denotes the space of bounded linear operators from $E_{1}$ into $E_{2}.$ For 
$E_{1}=E_{2}=E$ it will be denoted by $L\left( E\right) .$

A linear operator\ $A$ is said to be $\varphi $-positive in a Banach\ space $%
E$ with bound $M>0$ if $D\left( A\right) $ is dense on $E$ and $\left\Vert
\left( A+\lambda I\right) ^{-1}\right\Vert _{L\left( E\right) }\leq M\left(
1+\left\vert \lambda \right\vert \right) ^{-1}$ $\ $for any $\lambda \in
S_{\varphi },$ $0\leq \varphi <\pi ,$ where $I$ is the identity operator in $%
E$. Sometimes $A+\lambda I$\ will be written as $A+\lambda $ and will be
denoted by $A_{\lambda }.$ It is known $\left[ \text{24, \S 1.15.1}\right] $
that a positive operator $A$ has well-defined fractional powers\ $A^{\theta
}.$ Let $E\left( A^{\theta }\right) $ denote the space $D\left( A^{\theta
}\right) $ with norm 
\begin{equation*}
\left\Vert u\right\Vert _{E\left( A^{\theta }\right) }=\left( \left\Vert
u\right\Vert ^{p}+\left\Vert A^{\theta }u\right\Vert ^{p}\right) ^{\frac{1}{p%
}},1\leq p<\infty ,\text{ }0<\theta <\infty .
\end{equation*}

Let $S\left( R^{n};E\right) $ denote the Schwartz class, i.e., the space of
all $E$-valued rapidly decreasing smooth functions on $R^{n}$ and $C\left(
\Omega ;E\right) $ denotes the space of all $E$-valued norm bounded
functions on $\Omega .$ Let $F$ denote the Fourier transformation. A
function $\Psi \in C\left( R^{n};L\left( E\right) \right) $ is called
Fourier multiplier in $L_{p,\gamma }\left( R^{n};E\right) $ if the map 
\begin{equation*}
u\rightarrow \Phi u=F^{-1}\Psi \left( \xi \right) Fu,\text{ }u\in S\left(
R^{n};E\right)
\end{equation*}%
is well defined and extends to a bounded linear operator in $L_{p}\left(
R^{n};E\right) .$ The set of all multipliers in\ $L_{p}\left( R^{n};E\right) 
$\ will denoted by $M_{p}^{p}\left( E\right) .$

\textbf{Definition 1}$.$\textbf{1}$.$ A Banach space $E$\ is said to be a
space satisfying multiplier condition with respect to $p\in \left( 1,\infty
\right) $ if for any $\Psi \in C^{\left( 1\right) }\left( \mathbb{R};L\left(
E\right) \right) $ the $R$-boundedness (see e.g. $\left[ \text{9, \S\ 4.1}%
\right] $) of the set%
\begin{equation*}
\left\{ \xi ^{j}\Psi ^{\left( j\right) }\left( \xi \right) :\xi \in \mathbb{R%
}\backslash \left\{ 0\right\} ,\text{ }j=0,1\right\}
\end{equation*}%
implies $\Psi \in $ $M_{p}^{p}\left( E\right) .$

\textbf{Remark 1.1.\ }Note that if $E$ is $UMD$ space, then for example, by $%
\left[ \text{25}\right] $, $\left[ \text{9}\right] $, $\left[ \text{11}%
\right] $ this space satisfies the multiplier condition.

By $\left( E_{1},E_{2}\right) _{\theta ,p}$, $0<\theta <1,1\leq p\leq \infty 
$ we will denote the interpolation spaces obtained from $\left\{
E_{1},E_{2}\right\} $ by the $K$-method \ $\left[ \text{24, \S 1.3.2}\right] 
$.\ 

The operator $A\left( x\right) $ is said to be $\varphi $-positive uniformly
with respect to $x\in G$ in\ $E$ with bound $M>0$ if $D\left( A\left(
x\right) \right) $ is independent of$\ x$, $D\left( A\left( x\right) \right) 
$ is dense in $E$ and $\left\Vert \left( A\left( x\right) +\lambda \right)
^{-1}\right\Vert \leq \frac{M}{1+\left\vert \lambda \right\vert }$ for all $%
\lambda \in S\left( \varphi \right) ,$ $0\leq \varphi <\pi $, where $M$ is
independent of $x.$

The $\varphi $-positive operator $A\left( x\right) ,$ $x\in \sigma $ is said
to be uniformly $R$-positive in a Banach space $E$ if there exists $\varphi
\in \left[ 0\right. ,\left. \pi \right) $ such that the set 
\begin{equation*}
\left\{ A\left( x\right) \left( A\left( x\right) +\xi I\right) ^{-1}:\xi \in
S_{\varphi }\right\}
\end{equation*}%
is uniformly $R$-bounded, that is 
\begin{equation*}
\text{ }\sup\limits_{x\in \sigma }R\left( \left\{ \left[ A\left( x\right)
\left( A\left( x\right) +\xi I\right) ^{-1}\right] \text{: }\xi \in
S_{\varphi }\right\} \right) \leq M.
\end{equation*}

\bigskip\ Let $E_{0}$ and $E$ be two Banach spaces and $E_{0}$ is
continuously and densely embedded into $E$. Let $\sigma $ be a domi\i n in $%
\mathbb{R}.$ Consider the Sobolev-Lions type space\ $W_{p}^{m}\left( \sigma
;E_{0},E\right) $ that consisting of all functions $u\in L_{p}\left( \sigma
;E_{0}\right) $ that have generalized derivatives $u^{\left( m\right) }\in
L_{p}\left( \sigma ;E\right) $ with the norm 
\begin{equation*}
\ \left\Vert u\right\Vert _{W_{p}^{m}}=\left\Vert u\right\Vert
_{W_{p}^{m}\left( \sigma ;E_{0},E\right) }=\left\Vert u\right\Vert
_{L_{p}\left( \sigma ;E_{0}\right) }+\left\Vert u^{\left( m\right)
}\right\Vert _{L_{p}\left( \sigma ;E\right) }<\infty .
\end{equation*}

The embedding theorems play a key role in the perturbation theory of DOEs.
For estimating lower order derivatives we use following embedding theorems
from $\left[ 21\right] $:

\textbf{Theorem A}$_{1}$. Assume the following conditions are satisfied:

(1) $E$ is a Banach space satisfying the multiplier condition with respect
to $p;$

(2)\ $A$ is an $R$-positive operator in $E,$ $\sigma \subset \mathbb{R};$

(3)\ $0\leq j\leq m,$ $0\leq \mu \leq 1-\frac{j}{m}$, $1<p<\infty $; $h$ is
a positive parameter that $0<h<h_{0}<\infty ;$

(4)\ There exists a bounded linear extension operator from $W_{p}^{m}\left(
\sigma ;E\left( A\right) ,E\right) $ to $W_{p}^{m}\left( \left( -\infty
,\infty \right) ;E\left( A\right) ,E\right) $.

Then the embedding $D^{j}W_{p}^{m}\left( \sigma ;E\left( A\right) ,E\right)
\subset L_{p}\left( \sigma ;E\left( A^{1-\frac{j}{m}-\mu }\right) \right) $
is continuous. Moreover, for $u\in W_{p}^{m}\left( \sigma ;E\left( A\right)
,E\right) $ the following estimate holds 
\begin{equation*}
\left\Vert u^{\left( j\right) }\right\Vert _{L_{p}\left( \sigma ;E\left(
A^{1-\frac{j}{m}-\mu }\right) \right) }\leq h^{\mu }\left\Vert u\right\Vert
_{W_{p}^{m}\left( \sigma ;E\left( A\right) ,E\right) }+h^{-\left( 1-\mu
\right) }\left\Vert u\right\Vert _{L_{p}\left( \sigma ;E\right) }.
\end{equation*}

Consider the DOE with variable coefficients on $\left( -\infty ,\infty
\right) $

\begin{equation}
\ \varepsilon a\left( x\right) u^{\left( 2\right) }\left( x\right) +A\left(
x\right) u\left( x\right) +\dsum\limits_{i=0}^{1}\varepsilon ^{\frac{i}{2}%
}A_{i}\left( x\right) u^{\left( i\right) }\left( x\right) +\lambda u\left(
x\right) =f\left( x\right) ,\text{ }  \tag{1.6}
\end{equation}%
where $a\left( .\right) $ is a real-valued function, $\varepsilon $ is a
positive parameter, $A\left( .\right) $ and \ $A_{j}\left( .\right) $ are
linear operator functions in a Banach space $E,$ $\lambda $ is a complex
parameter$.$

Let $\omega _{1}=\omega _{1}\left( x\right) $, $\omega _{2}=\omega
_{2}\left( x\right) $ be roots of the equation $a\left( x\right) \omega
^{2}+1=0$.

From $\left[ \text{21}\right] $ we obtain

\textbf{Theorem A}$_{2}$\textbf{. }Suppose the following conditions are
satisfied:

(1) $E$ is a Banach space satisfying the multiplier condition with respect to%
$\ p\in \left( 1,\infty \right) ;$

(2) $A\left( x\right) $ is an $R$-positive operator in $E$ for $\varphi \in %
\left[ 0,\right. \left. \pi \right) $ uniformly with respect to $x\in \left[
0,1\right] $ and $A\left( x\right) A^{-1}\left( x_{0}\right) \in C\left(
\left( -\infty ,\infty \right) ;L\left( E\right) \right) $ for a.e. $%
x_{0}\in \left( -\infty ,\infty \right) ;$

(3) for any $\delta >0$ there is a positive $C\left( \delta \right) $ such
that%
\begin{equation*}
\left\Vert A_{1}\left( x\right) u\right\Vert \leq \delta \left\Vert
u\right\Vert _{\left( E\left( A\right) ,E\right) _{\frac{1}{2},\infty
}}+C\left( \delta \right) \left\Vert u\right\Vert
\end{equation*}%
for $u\in \left( E\left( A\right) ,E\right) _{\frac{1}{2},\infty }$ and $%
\left\Vert A_{0}\left( x\right) u\right\Vert \leq \delta \left\Vert
Au\right\Vert _{E}+C\left( \delta \right) \left\Vert u\right\Vert $ for $%
u\in D\left( A\right) ;$

(4) $a\in C_{b}\left( -\infty ,\infty \right) $ and $\func{Re}\omega
_{k}\neq 0$ and $\frac{\lambda }{\omega _{k}}\in S\left( \varphi \right) $
for $\lambda \in S\left( \varphi \right) $, $0\leq \varphi <\pi ,$ $k=1,2$.$%
\ $a.e. $x\in \mathbb{R};$

Then problem $\left( 1.6\right) $ has a unique solution $u\in
W_{p}^{2}\left( \mathbb{R};E\left( A\right) ,E\right) $ for\ $f\in
L_{p}\left( \mathbb{R};E\right) .$ Moreover, for $\left\vert \arg \lambda
\right\vert \leq \varphi $ and sufficiently large $\left\vert \lambda
\right\vert $ the following uniform coercive estimate holds

\begin{equation*}
\sum\limits_{i=0}^{2}\left\vert \lambda \right\vert ^{1-\frac{i}{2}%
}\varepsilon ^{\frac{i}{2}}\left\Vert u^{\left( i\right) }\right\Vert
_{L_{p}\left( \mathbb{R};E\right) }+\left\Vert Au\right\Vert _{L_{p}\left( 
\mathbb{R};E\right) }\leq C\left\Vert f\right\Vert _{L_{p}\left( \mathbb{R}%
;E\right) }.
\end{equation*}

Consider the nonhomogenous BVP for DOE with constant coefficients on half
plane

\begin{equation}
\ \varepsilon au^{\left( 2\right) }\left( x\right) +Au\left( x\right)
+\lambda u\left( x\right) =f\left( x\right) ,\text{ }x\in \left( 0,\infty
\right) ,  \tag{1.7}
\end{equation}%
\begin{equation*}
\dsum\limits_{i=0}^{\nu }\alpha _{i}\varepsilon ^{\nu _{i}}u^{\left(
i\right) }\left( 0\right) =\varkappa ,
\end{equation*}%
where $\ \varkappa \in \left( E\left( A\right) ,E\right) _{\frac{\nu }{2},p}$%
, $a$ is a complex number, $\varepsilon $ is a positive parameter, $\nu _{i}=%
\frac{i}{2}+\frac{1}{2p};$ $A$ is a linear operator in a Banach space $E,$ $%
\lambda $ is a complex parameter, $\alpha _{i}$ are complex numbers and $\nu
\in \left\{ 0,1\right\} $, $\alpha _{\nu }\neq 0.$

Let $\omega _{1}$, $\omega _{2}$ be roots of equation $a\omega ^{2}+1=0$.

From $\left[ 22\right] $ we obtain.

\textbf{Theorem A}$_{3}$\textbf{. }Suppose the following conditions are
satisfied:

(1) $E$ is a Banach space satisfying the multiplier condition with respect to%
$\ p\in \left( 1,\infty \right) ;$

(2) $A$ is an $R$-positive operator in $E$ for $\varphi \in \left[ 0,\right.
\left. \pi \right) ;$

(4) $a$ is a complex number such that $\func{Re}\omega _{k}\neq 0$ and $%
\frac{\lambda }{\omega _{k}}\in S\left( \varphi \right) $ for $\lambda \in
S\left( \varphi \right) $, $0\leq \varphi <\pi ,$ $k=1,2$.

Then problem $\left( 1.7\right) $ has a unique solution $u\in
W_{p}^{2}\left( 0,\infty ;E\left( A\right) ,E\right) $ for\ $f\in
L_{p}\left( 0,\infty ;E\right) .$ Moreover, for $\left\vert \arg \lambda
\right\vert \leq \varphi $ and sufficiently large $\left\vert \lambda
\right\vert $ the following uniform coercive estimate holds

\begin{equation*}
\sum\limits_{i=0}^{2}\left\vert \lambda \right\vert ^{1-\frac{i}{2}%
}\varepsilon ^{\frac{i}{2}}\left\Vert u^{\left( i\right) }\right\Vert
_{L_{p}\left( 0,\infty ;E\right) }+\left\Vert Au\right\Vert _{L_{p}\left(
0,\infty ;E\right) }\leq C\left\Vert f\right\Vert _{L_{p}\left( 0,\infty
;E\right) }+\left\Vert \varkappa \right\Vert _{\left( E\left( A\right)
,E\right) _{\frac{\nu }{2},p}}.
\end{equation*}

Consider the nonlocal BVP for DOE with constant coefficients

\begin{equation}
\ \varepsilon au^{\left( 2\right) }\left( x\right) +Au\left( x\right)
+\lambda u\left( x\right) =f\left( x\right) ,\text{ }  \tag{1.8}
\end{equation}%
\begin{equation*}
\sum\limits_{i=0}^{\mu _{k}}\varepsilon ^{\nu _{i}}\left[ \alpha
_{ki}u^{\left( i\right) }\left( 0\right) +\beta _{ki}u^{\left( i\right)
}\left( 1\right) \right] =f_{k},\text{ }k=1,2,\text{ }x\in \left( 0,1\right)
,
\end{equation*}%
where $\ \ f_{k}\in \left( E\left( A\right) ,E\right) _{\frac{p\mu _{k}+1}{2p%
},p}$, $A$ is a linear operator in a Banach space $E,$ $\varepsilon $ is a
positive parameter, $\nu _{i}=\frac{i}{2}+\frac{1}{2p}$, $\lambda $ is a
complex parameter, $a,$ $\alpha _{ki},$ $\beta _{ki}$ are complex numbers
and $\mu _{k}\in \left\{ 0,1\right\} .$

From $\left[ \text{20}\right] $ we obtain.

\textbf{Theorem A}$_{4}.$ Suppose the following conditions are satisfied:

(1) $E$ is a Banach space satisfying the multiplier condition with respect to%
$\ p\in \left( 1,\infty \right) ;$

(2) $A$ is an $R$-positive operator in $E$ for $\varphi \in \left[ 0,\right.
\left. \pi \right) ;$

(3) $a$ is a complex number such that $\func{Re}\omega _{k}\neq 0$ and $%
\frac{\lambda }{\omega _{k}}\in S\left( \varphi \right) $ for $\lambda \in
S\left( \varphi \right) $, $0\leq \varphi <\pi ,$ $k=1,2$;

(4) $\alpha _{k}=\alpha _{k\nu _{k}}\neq 0$, $\beta _{k}=\beta _{k\nu
_{k}}\neq 0$, $\eta =\left( -1\right) ^{\mu _{1}}\alpha _{1}\beta
_{2}-\left( -1\right) ^{\mu _{2}}\alpha _{2}\beta _{1}\neq 0$, $a>0;$

Then problem $\left( 1.8\right) $ has a unique solution $u\in
W_{p}^{2}\left( 0,1;E\left( A\right) ,E\right) $ for\ $f\in L_{p}\left(
0,1;E\right) $ and $f_{k}\in \left( E\left( A\right) ,E\right) _{\frac{p\nu
_{k}+1}{2},p}.$ Moreover, for $\left\vert \arg \lambda \right\vert \leq
\varphi $ and sufficiently large $\left\vert \lambda \right\vert $ the
following uniform coercive estimate holds

\begin{equation*}
\sum\limits_{i=0}^{2}\left\vert \lambda \right\vert ^{1-\frac{i}{2}%
}\varepsilon ^{\frac{i}{2}}\left\Vert u^{\left( i\right) }\right\Vert
_{L_{p}\left( 0,1;E\right) }+\left\Vert Au\right\Vert _{L_{p}\left(
0,1;E\right) }\leq
\end{equation*}

\begin{equation*}
C\left[ \left\Vert f\right\Vert _{L_{p}\left( 0,1;E\right)
}+\sum\limits_{k=1}^{2}\left\Vert f_{k}\right\Vert _{\left( E\left( A\right)
,E\right) _{\frac{p\mu _{k}+1}{2},p}}\right] .
\end{equation*}

By virtue Lions Petree trace theorem (see of $\left[ \text{24, \S 1.8.2}%
\right] $) we obtain

\textbf{Theorem A}$_{5}.$ Assume $m$ and $j$ are integers, $0\leq j\leq m-1,$
$\theta _{j}=\frac{pj+1}{pm}$, $p\in \left( 1,\infty \right) ;$ $\varepsilon
\in \left( 0,1\right) $ is a parameter$,$ $x_{0}\in \left[ 0,b\right] $.
Then, the linear transformation $u\rightarrow u^{\left( j\right) }\left(
x_{0}\right) $ is bounded from $W_{p}^{m}\left( 0,b;E_{0},E\right) $ onto $%
\left( E_{0},E\right) _{\theta _{j},p}$ and the following inequality holds \ 
\begin{equation*}
\varepsilon ^{\theta _{j}}\left\Vert u^{\left( j\right) }\left( x_{0}\right)
\right\Vert _{\left( E_{0},E\right) _{\theta _{j},p}}\leq C\left( \left\Vert
\varepsilon u^{\left( m\right) }\right\Vert _{L_{p,\gamma }\left(
0,b;E\right) }+\left\Vert u\right\Vert _{L_{p,\gamma \ }\left(
0,b;E_{0}\right) }\right) .
\end{equation*}

\begin{center}
\textbf{2. Abstract equation with variable coefficients }
\end{center}

\ \ Consider the exterior BVP for differential-operator equation with
variable coefficients%
\begin{equation}
Lu=\varepsilon au^{\left( 2\right) }+Au+\dsum\limits_{i=0}^{1}\varepsilon ^{%
\frac{i}{2}}A_{i}u^{\left( i\right) }+\lambda u=f,  \tag{2.1}
\end{equation}%
\begin{equation}
L_{1}u=\sum\limits_{i=0}^{\mu _{1}}\alpha _{i}\varepsilon ^{\nu
_{i}}u^{\left( i\right) }\left( 0\right) =0,\text{ }L_{2}u=\sum%
\limits_{i=0}^{\mu _{2}}\beta _{i}\varepsilon ^{\nu _{i}}u^{\left( i\right)
}\left( b\right) =0\text{, }  \tag{2.2}
\end{equation}%
where $a=a\left( x\right) $ is a complex-valued function, $\varepsilon $ is
a positive parameter, $\nu _{i}=\frac{i}{2}+\frac{1}{2p}$, $u=u\left(
x\right) $, $f=f\left( x\right) ,$ $x\in \sigma $ are $E$-valued unknown and
date functions; $A=A\left( x\right) $ and \ $A_{j}=A_{j}\left( x\right) $
are linear operator functions in a Banach space $E,$ $\lambda $ is a complex
parameter, $\alpha _{i},$ $\beta _{i}$ are complex numbers, $\mu _{k}\in
\left\{ 0,1\right\} $ and $\sigma =\mathbb{R}\setminus \left[ 0,b\right] .$

A function \ $u\in $ $W_{p}^{2}\left( \sigma ;E\left( A\right) ,E\right) $\
satisfying the equation $\left( 2.1\right) $ a.e. on $\sigma $ is said to be
the solution of the equation $\left( 2.1\right) $ on $\sigma .$

Consider the problem $\left( 2.1\right) -\left( 2.2\right) .$ Let $%
X=L_{p}\left( \sigma ;E\right) $ and $Y=W_{p}^{2}\left( \sigma ;E\left(
A\right) ,E\right) .$ Let $\omega _{1}=\omega _{1}\left( x\right) $, $\omega
_{2}=\omega _{2}\left( x\right) $ be roots of equation $a\left( x\right)
\omega ^{2}+1=0$.

The main result of this section is the following:

\textbf{Theorem 2.1. }Assume the following conditions are satisfied:

Suppose the following conditions are satisfied:

(1) $E$ is a Banach space satisfying the multiplier condition with respect to%
$\ p\in \left( 1,\infty \right) ;$

(2) $A\left( x\right) $ is an $R$-positive operator in $E$ for $\varphi \in %
\left[ 0,\right. \left. \pi \right) $ uniformly with respect to $x\in \left[
0,1\right] $ and $A\left( x\right) A^{-1}\left( x_{0}\right) \in C\left( 
\bar{\sigma};L\left( E\right) \right) $ for $x_{0}\in \left( 0,1\right) ;$

(3) for any $\delta >0$ there is a positive $C\left( \delta \right) $ such
that

$\left\Vert A_{1}\left( x\right) u\right\Vert \leq \delta \left\Vert
u\right\Vert _{\left( E\left( A\right) ,E\right) _{\frac{1}{2},\infty
}}+C\left( \delta \right) \left\Vert u\right\Vert $ for $u\in \left( E\left(
A\right) ,E\right) _{\frac{1}{2},\infty }$ and 
\begin{equation*}
\left\Vert A_{0}\left( x\right) u\right\Vert \leq \delta \left\Vert
Au\right\Vert _{E}+C\left( \delta \right) \left\Vert u\right\Vert
\end{equation*}%
for $u\in D\left( A\right) ;$

(4) $a\in C\left( \bar{\sigma}\right) $, $\func{Re}\omega _{k}\neq 0$ and $%
\frac{\lambda }{\omega _{k}}\in S\left( \varphi \right) $ for $\lambda \in
S\left( \varphi \right) $, $0\leq \varphi <\pi ,$ $k=1,2$.$\ $a.e. $x\in 
\mathbb{\sigma }.$

Then problem $\left( 2.1\right) -\left( 2.2\right) $ has a unique solution $%
u\in W_{p}^{2}\left( \sigma ;E\left( A\right) ,E\right) $ for\ $f\in
L_{p}\left( \sigma ;E\right) .$ Moreover, for $\left\vert \arg \lambda
\right\vert \leq \varphi $ and sufficiently large $\left\vert \lambda
\right\vert $ the following uniform coercive estimate holds

\begin{equation}
\sum\limits_{i=0}^{2}\left\vert \lambda \right\vert ^{1-\frac{i}{2}%
}\varepsilon ^{\frac{i}{2}}\left\Vert u^{\left( i\right) }\right\Vert
_{L_{p}\left( \sigma ;E\right) }+\left\Vert Au\right\Vert _{L_{p}\left(
\sigma ;E\right) }\leq C\left\Vert f\right\Vert _{L_{p}\left( \sigma
;E\right) }.  \tag{2.3}
\end{equation}

\textbf{Proof.} First of all, we will show the uniqueness of solution. Let $%
G_{1,}G_{2},...,G_{n}...$ be regions in $\mathbb{R}$ and $\varphi
_{1},\varphi _{2},...,\varphi _{n}...$ correspond to a partition of unit on $%
\sigma $, which functions $\varphi _{j}$ are smooth functions on $\mathbb{R}$%
, supp $\varphi _{j}\subset G_{j}$ and $\sum\limits_{j=1}^{\infty }\varphi
_{j}\left( x\right) =1$ for $x\in \sigma .$ Then for all $u\in Y$ we have $%
u\left( x\right) =$ $\sum\limits_{j=1}^{\infty }u_{j}\left( x\right) ,$
where $u_{j}\left( x\right) =u\left( x\right) \varphi \left( x\right) .$ Let 
$u\in Y$ be a solution of $\left( 2.1\right) -\left( 2.2\right) .$ Then from 
$\left( 2.1\right) -\left( 2.2\right) $ we obtain

\begin{equation}
\left( \ L+\lambda \right) u_{j}=\varepsilon au_{j}^{\left( 2\right) }\left(
x\right) +\left( A+\lambda \right) u_{j}\left( x\right) =f_{j}\left(
x\right) ,  \tag{2.4}
\end{equation}%
where

\begin{equation}
f_{j}=f\varphi _{j}+\varepsilon a\left( 2u^{\left( 1\right) }\varphi
_{j}^{\left( 1\right) }+u\varphi _{j}^{\left( 2\right) }\right) +\varepsilon
^{\frac{1}{2}}\varphi _{j}^{\left( 1\right) }A_{1}u,  \tag{2.5}
\end{equation}

\begin{equation*}
L_{1}u_{j}=\varkappa _{1},\text{ }L_{2}u_{j}=\varkappa _{2},\text{ }%
j=1,2,...,\infty ,
\end{equation*}%
\begin{equation*}
\varkappa _{1}=\alpha _{1}u\left( 0\right) \varphi _{j}\left( 0\right) ,%
\text{ }\varkappa _{2}=\alpha _{1}u\left( b\right) \varphi _{j}\left(
b\right) .
\end{equation*}%
By Lemma A$_{5},$ $\varkappa _{1},$ $\varkappa _{2}\in \left( E\left(
A\right) ,E\right) _{\frac{1}{2p},p}$. By freezing coefficients in $\left(
2.4\right) $ we obtain that

\begin{equation}
\bigskip \varepsilon a\left( x_{0j}\right) u_{j}^{\left( 2\right) }\left(
x\right) +\left( A\left( x_{0j}\right) +\lambda \right) u_{j}\left( x\right)
=F_{j}\left( x\right) ,  \tag{2.6}
\end{equation}%
\begin{equation*}
L_{1}u_{j}=\varkappa _{1},\text{ }L_{2}u_{j}=\varkappa _{2},\text{ }%
j=1,2,...,\infty ,
\end{equation*}%
where 
\begin{equation}
F_{j}=f_{j}+\left[ A\left( x_{0j}\right) -A\left( x\right) \right]
u_{j}-\varepsilon \left[ a\left( x\right) -a\left( x_{0j}\right) \right]
u_{j}^{\left( 2\right) }.  \tag{2.7}
\end{equation}%
Since functions $u_{j}\left( x\right) $ have compact supports, by extending $%
u_{j}\left( x\right) $ on the outsides of supp $\varphi _{j}$ we obtain BVPs
for DOEs with constant coefficients

\begin{equation}
\varepsilon a\left( x_{0j}\right) u_{j}^{\left( 2\right) }+\left( A\left(
x_{0j}\right) +\lambda \right) u_{j}=F_{j},\text{ }L_{1}u_{j}=\varkappa _{1},%
\text{ }L_{2}u_{j}=\varkappa _{2}.  \tag{2.8}
\end{equation}%
Since $a$ is uniformly bounded on $\sigma $ for all small $\rho >0$ there is
a large $r_{0}>0$ such that $\left\vert a\left( x\right) -a\left( \pm \infty
\right) \right\vert \leq \delta $ for all $\left\vert x\right\vert \geq
r_{0}.$ Let 
\begin{equation*}
G_{0}=\left( -\infty ,\infty \right) \setminus O_{r_{0}}\left( 0\right) 
\text{, }O_{r_{0}}\left( 0\right) =\left\{ x\in \sigma ,\text{ }\left\vert
x\right\vert \leq r_{0}\right\} .
\end{equation*}%
Cover $O_{r_{0}}\left( 0\right) $ by finitely many intervals $%
G_{j}=O_{r_{j}}\left( x_{0j}\right) $ such that%
\begin{equation*}
\left\vert a\left( x\right) -a\left( x_{0j}\right) \right\vert \leq \delta 
\text{ for }\left\vert x-x_{0j}\right\vert \leq r_{j},j=1,2,....
\end{equation*}%
Define coefficients of local operators, i.e. 
\begin{eqnarray*}
a^{0}\left( x\right) &=&\left\{ 
\begin{array}{c}
a\left( x\right) \text{, }x\notin O_{r_{0}}\left( 0\right) \\ 
a\left( r_{0}^{2}\frac{x}{\left\vert x\right\vert ^{2}}\right) ,\text{ }x\in
O_{r_{0}}\left( 0\right) \text{\ \ \ \ \ \ \ }%
\end{array}%
\right\} , \\
a^{j}\left( x\right) &=&\left\{ 
\begin{array}{c}
a\left( x\right) \text{, }x\in O_{r_{j}}\left( x_{0j}\right) \\ 
a\left( x_{0j}+r_{0}^{2}\frac{x-x_{0j}}{\left\vert x-x_{0j}\right\vert ^{2}}%
\right) \text{, }x\notin O_{r_{j}}\left( x_{0j}\right) \text{\ \ \ \ \ \ \ }%
\end{array}%
\right\}
\end{eqnarray*}%
and 
\begin{eqnarray*}
A^{0}\left( x\right) A^{-1}\left( x_{0j}\right) &=&\left\{ 
\begin{array}{c}
A\left( x\right) A^{-1}\left( x_{0j}\right) \text{, }x\notin O_{r_{0}}\left(
0\right) \\ 
A\left( r_{0}^{2}\frac{x}{\left\vert x\right\vert ^{2}}\right) A^{-1}\left(
x_{0j}\right) ,\text{ }x\in O_{r_{0}}\left( 0\right) \text{\ \ \ \ \ \ \ }%
\end{array}%
\right\} , \\
A^{j}\left( x\right) A^{-1}\left( x_{0j}\right) &=&\left\{ 
\begin{array}{c}
A\left( x\right) A^{-1}\left( x_{0j}\right) \text{, }x\in O_{r_{j}}\left(
x_{0j}\right) \\ 
A\left( x_{0j}+r_{0}^{2}\frac{x-x_{0j}}{\left\vert x-x_{0j}\right\vert ^{2}}%
\right) A^{-1}\left( x_{0j}\right) \text{, }y\notin O_{r_{j}}\left(
x_{0j}\right) \text{\ \ \ \ \ \ \ }%
\end{array}%
\right\}
\end{eqnarray*}%
for each $j=1,2,....$ Then, for all $x\in \sigma $ and $j=0,1,2,....$we get%
\begin{equation*}
\left\vert a^{j}\left( x\right) -a^{j}\left( x_{0j}\right) \right\vert \leq
\delta \text{ and }\left\Vert A^{j}\left( x\right) A^{-1}\left(
x_{0j}\right) -A^{j}\left( x_{0j}\right) A^{-1}\left( x_{0j}\right)
\right\Vert _{B\left( E\right) }<\delta .
\end{equation*}%
Let $\varphi _{j}$ such that $0$, $b\in $ supp $\varphi _{j}.$ Then by
virtue of Theorem A$_{4}$ we obtain that problem $\left( 2.8\right) $ has a
unique solution $u_{j}$\ \bigskip and the coercive uniform estimates hold 
\begin{equation*}
\sum\limits_{i=0}^{2}\left\vert \lambda \right\vert ^{1-\frac{i}{2}%
}\varepsilon ^{\frac{i}{2}}\left\Vert u_{j}^{\left( i\right) }\right\Vert
_{G_{j},p}+\left\Vert Au_{j}\right\Vert _{G_{j},p}\leq C\left\Vert
F_{j}\right\Vert _{G_{j},p}+\dsum\limits_{k=1}^{2}\left\Vert \varkappa
_{k}\right\Vert _{E_{p}},
\end{equation*}%
where, $\left\Vert .\right\Vert _{G_{j},p}$ denote $E$-valued $L_{p}$-norms
on $G_{j}$ and $E_{p}=\left( E\left( A\right) ,E\right) _{\frac{1}{2p},p}.$
Then by using Theorems A$_{1}$ and A$_{6}$ we obtain from the above estimate
the following%
\begin{equation}
\sum\limits_{i=0}^{2}\left\vert \lambda \right\vert ^{1-\frac{i}{2}%
}\varepsilon ^{\frac{i}{2}}\left\Vert u_{j}^{\left( i\right) }\right\Vert
_{G_{j},p}+\left\Vert Au_{j}\right\Vert _{G_{j},p}\leq C\left\Vert
F_{j}\right\Vert _{G_{j},p}.  \tag{2.9}
\end{equation}%
Let $\varphi _{j}$ such that $0$, $1\bar{\in}$ supp $\varphi _{j}.$ Hence, $%
\varkappa _{k}=0.$ Then in a similar way, Theorem A$_{2}$ and Theorem A$_{3}$
imply the same estimates 
\begin{equation}
\sum\limits_{i=0}^{2}\left\vert \lambda \right\vert ^{1-\frac{i}{2}%
}\varepsilon ^{\frac{i}{2}}\left\Vert u_{j}^{\left( i\right) }\right\Vert
_{G_{j},p}+\left\Vert Au_{j}\right\Vert _{G_{j},p}\leq C\left\Vert
F_{j}\right\Vert _{G_{j},p}  \tag{2.10}
\end{equation}%
for domains $G_{j}$ \ adjoin the boundary point $0$ and $b$. Hence, using
properties of the smoothness of coefficients of equations $\left( 2.5\right)
,\left( 2.7\right) $ and choosing diameters of supp$\varphi _{j}$
sufficiently small, we get 
\begin{equation}
\left\Vert F_{j}\right\Vert _{G_{j},p}\leq \delta \left\Vert
u_{j}\right\Vert _{W_{p}^{2}\left( G_{j};E\left( A\right) ,E\right)
}+C\left( \delta \right) \left\Vert f_{j}\right\Vert _{G_{j},p},\text{ } 
\tag{2.11}
\end{equation}%
where $\delta $ is a sufficiently small positive number and $C\left( \delta
\right) $ is a continuous function.\ Consequently, from $\left( 2.9\right) $-%
$\left( 2.11\right) $ by using Theorem A$_{1}$ we get 
\begin{equation*}
\sum\limits_{i=0}^{2}\left\vert \lambda \right\vert ^{1-\frac{i}{2}%
}\varepsilon ^{\frac{i}{2}}\left\Vert u_{j}^{\left( i\right) }\right\Vert
_{G_{j},p}+\left\Vert Au_{j}\right\Vert _{G_{j},p}\leq
\end{equation*}

\begin{equation*}
C\left\Vert f\right\Vert _{G_{j},p}+\delta \left\Vert u_{j}\right\Vert
_{W_{p,\gamma }^{2}\left( G_{j};E\left( A\right) ,E\right) }+C\left( \delta
\right) \left\Vert u_{j}\right\Vert _{G_{j},p}.
\end{equation*}%
Choosing $\delta <1$ from the above inequality we have

\begin{equation}
\sum\limits_{i=0}^{2}\left\vert \lambda \right\vert ^{1-\frac{i}{2}%
}\varepsilon ^{\frac{i}{2}}\left\Vert u_{j}^{\left( i\right) }\right\Vert
_{G_{j},p}+\left\Vert Au_{j}\right\Vert _{G_{j},p}\leq C\left[ \left\Vert
f\right\Vert _{G_{j},p}+\left\Vert u_{j}\right\Vert _{G_{j},p}\right] \text{%
, }j=1,2,....  \tag{2.12}
\end{equation}%
Then using the equality $u\left( x\right) =$ $\sum\limits_{j=1}^{\infty
}u_{j}\left( x\right) $ and the estimate $\left( 2.12\right) $ for $u\in Y$
we have 
\begin{equation}
\bigskip \sum\limits_{i=0}^{2}\left\vert \lambda \right\vert ^{1-\frac{i}{2}%
}\varepsilon ^{\frac{i}{2}}\left\Vert u^{\left( i\right) }\right\Vert
_{p}+\left\Vert Au_{j}\right\Vert _{p}\leq C\left[ \left\Vert \left(
L+\lambda \right) u\right\Vert _{p}+\left\Vert u\right\Vert _{p}\right] . 
\tag{2.13}
\end{equation}%
Let $u\in Y$ be solution of problem $\left( 2.1\right) -\left( 2.2\right) .$
Then for $\left\vert \arg \lambda \right\vert \leq \varphi $ we have

\begin{equation}
\left\Vert u\right\Vert _{X}=\left\Vert \left( L+\lambda \right)
u-Lu\right\Vert _{X}\leq \frac{1}{\lambda }\left[ \left\Vert \left(
L+\lambda \right) u\right\Vert _{X}+\left\Vert u\right\Vert _{Y}\right] . 
\tag{2.14}
\end{equation}%
Then by Theorem A$_{1}$, by virtue of $\left( 2.12\right) $ and $\left(
2.14\right) $ for sufficiently large $\left\vert \lambda \right\vert $ we
have

\begin{equation}
\sum\limits_{i=0}^{2}\left\vert \lambda \right\vert ^{1-\frac{i}{2}%
}\varepsilon ^{\frac{i}{2}}\left\Vert u^{\left( i\right) }\right\Vert
_{X}+\left\Vert Au\right\Vert _{X}\leq C\left\Vert \left( L+\lambda \right)
u\right\Vert _{X}.  \tag{2.15}
\end{equation}%
Consider the operator $O_{\varepsilon }$ in $X$ generated by problem $\left(
2.1\right) -\left( 2.2\right) $, i.e., 
\begin{equation*}
D\left( O_{\varepsilon }\right) =W_{p}^{2}\left( \sigma ;E\left( A\right)
,E,L_{1},L_{2}\right) ,\text{ }
\end{equation*}%
\begin{equation*}
O_{\varepsilon }u=\varepsilon au^{\left( 2\right) }\left( x\right) +Au\left(
x\right) +\dsum\limits_{i=0}^{1}\varepsilon ^{\frac{i}{2}}A_{i}u^{\left(
i\right) }.
\end{equation*}%
The estimate $\left( 2.15\right) $ implies that the problem $\left(
2.1\right) -\left( 2.2\right) $ has only a unique solution and the operator $%
O_{\text{ }}+\lambda $ has an invertible operator in its rank space. We need
to show that this rank space coincides with the space $X.$ We consider the
smooth functions $g_{j}=g_{j}\left( x\right) $ with respect to the partition
of the unique $\varphi _{j}=\varphi _{j}\left( x\right) $ on $\sigma $ that
equal one on supp $\varphi _{j},$ where supp $g_{j}\subset G_{j}$ and $%
\left\vert g_{j}\left( x\right) \right\vert <1.$ Let us construct the
function $u_{j}$ for all $j$, that are defined on $\Omega _{j}=\sigma \cap
G_{j}$ and satisfying the problem $\left( 2.1\right) -\left( 2.2\right) .$
The problem $\left( 2.1\right) -\left( 2.2\right) $ can be expressed as\ 
\begin{equation}
\varepsilon a\left( x_{0j}\right) u_{j}^{\left( 2\right) }+\left( A\left(
x_{0j}\right) +\lambda \right) u_{j}=g_{j}\left\{ F_{j}+\left[ A\left(
x_{0j}\right) -A\left( x\right) \right] u_{j}\right.  \tag{2.16}
\end{equation}%
\ 

\begin{equation*}
\left. -\left[ a\left( x\right) -a\left( x_{0j}\right) \right] u_{j}\right\}
,\text{ }L_{1}u_{j}=0,\text{ }L_{2}u_{j}=0\text{, }j=1,2,....
\end{equation*}%
Consider the $L_{p}\left( G_{j};E\right) -$realization of the above local
operators $O_{j\lambda \varepsilon }=O_{\varepsilon j}+\lambda $ defined as%
\begin{equation*}
D\left( O_{\varepsilon j\lambda }\right) =W_{p}^{2}\left( G_{j};E\left(
A\right) ,E,L_{1},L_{2}\right) ,\text{ }
\end{equation*}%
\begin{equation*}
O_{\varepsilon j\lambda }u=\varepsilon a\left( x_{0j}\right) u^{\left(
2\right) }\left( x\right) +\left( A\left( x_{0j}\right) +\lambda \right)
u\left( x\right) .
\end{equation*}%
By virtue of Theorem A$_{1},$ for $f\in L_{p}\left( G_{j};E\right) $, $%
\left\vert \arg \lambda \right\vert \leq \varphi $ and sufficiently large $%
\left\vert \lambda \right\vert $ we have

\begin{equation}
\sum\limits_{i=0}^{2}\left\vert \lambda \right\vert ^{1-\frac{i}{2}%
}\varepsilon ^{\frac{i}{2}}\left\Vert \frac{d^{i}}{dx^{i}}O_{j\lambda
}^{-1}{}f\right\Vert _{p}+\left\Vert AO_{j\lambda }^{-1}{}f\right\Vert
_{p}\leq C\left\Vert {}f\right\Vert _{p}.  \tag{2.17}
\end{equation}%
Extending $u_{j}$ zero on the outside of supp$\varphi _{j}$ and passing
substitutions $u_{j}=O_{\varepsilon j\lambda }^{-1}{}\upsilon _{j}$ in $%
\left( 2.17\right) ,$ obtain equations with respect to $\upsilon _{j}.$

\begin{equation}
\upsilon _{j}=K_{\varepsilon j\lambda }\upsilon _{j}+g_{j}f,j=1,2,...,N. 
\tag{2.18}
\end{equation}%
By virtue of Theorem A$_{1}$ and estimate $\left( 2.17\right) $, in view of
the smoothness of the coefficients of the expression $K_{j\lambda },$ for
sufficiently large $\left\vert \lambda \right\vert $ we have $\left\Vert
K_{j\lambda }\right\Vert <\delta ,$ where $\delta $ is sufficiently small.
Consequently, equations $\left( 2.18\right) $ have unique solutions $%
\upsilon _{j}=\left[ I-K_{\varepsilon j\lambda }{}\right] ^{-1}g_{j}f$ .\
Moreover, 
\begin{equation*}
\left\Vert \upsilon _{j}\right\Vert _{X}=\left\Vert \left[ I-K_{\varepsilon
j\lambda }{}\right] ^{-1}g_{j}f\right\Vert _{X}\leq \left\Vert f\right\Vert
_{X}.
\end{equation*}%
Whence, $\left[ I-K_{\varepsilon j\lambda }{}\right] ^{-1}g_{j}$ are bounded
linear operators from $X$ to $L_{p}\left( G_{j};E\right) .$ Thus, we obtain
that 
\begin{equation*}
u_{j}=U_{\varepsilon j\lambda }f=O_{\varepsilon j\lambda }^{-1}{}\left[
I-K_{\varepsilon j\lambda }{}\right] ^{-1}g_{j}f
\end{equation*}%
are solutions of $\left( 2.18\right) $. Consider the linear operator $\left(
U_{\varepsilon }+\lambda \right) $ in $X$ such that 
\begin{equation*}
\left( U+\lambda \right) f=\sum\limits_{j=1}^{\infty }\varphi _{j}\left(
y\right) U_{\varepsilon j\lambda }f.
\end{equation*}%
It is clear from the constructions $U_{\varepsilon j\lambda \text{ }}$and
the estimate $\left( 2.17\right) $ that operators $U_{\varepsilon j\lambda }$
are bounded linear from $X$ to $Y$ and

\begin{equation}
\sum\limits_{i=0}^{2}\left\vert \lambda \right\vert ^{1-\frac{i}{2}%
}\varepsilon ^{\frac{i}{2}}\left\Vert \frac{d^{i}}{d^{i}}U_{\varepsilon
j\lambda }^{-1}{}f\right\Vert _{X}+\left\Vert AU_{\varepsilon j\lambda
}^{-1}{}f\right\Vert _{X}\leq C\left\Vert f\right\Vert _{X}.  \tag{2.19}
\end{equation}%
Therefore, $\left( U+\lambda \right) $ is a bounded linear operator from $%
L_{p}$ to $L_{p}.$ Let $O$ denote the operator in $X$ generated by BVP $%
\left( 2.1\right) -\left( 2.2\right) .$ Then act of $\left( O+\lambda
\right) $ to $u=\sum\limits_{j=1}^{\infty }\varphi _{j}U_{\varepsilon
j\lambda }f$ \ gives $\left( O+\lambda \right) u=f+\sum\limits_{j=1}^{\infty
}\Phi _{\varepsilon j\lambda }f,$ where $\Phi _{\varepsilon j\lambda }$ are
a linear combination of $U_{\varepsilon j\lambda }$ and $\frac{d}{dx}%
U_{\varepsilon j\lambda }$. By virtue of embedding Theorem A$_{1}$, the
estimate $\left( 2.19\right) $ and from the expression $\Phi _{\varepsilon
j\lambda }$ we obtain that operators $\Phi _{j\lambda }$ are bounded linear
from $X$ to $L_{p}\left( G_{j};E\right) $ and $\left\Vert \Phi _{\varepsilon
j\lambda }\right\Vert <1$. Therefore, there exists a bounded linear\
invertible operator$\left( I+\sum\limits_{j=1}^{\infty }\Phi _{\varepsilon
j\lambda }\right) ^{-1}.$ So, we obtain that the BVP $\left( 2.1\right)
-\left( 2.2\right) $ for $\ f\in X$ has a unique solution 
\begin{equation}
u\left( x\right) =\left( O_{\varepsilon }+\lambda \right) ^{-1}f=\left(
U_{\varepsilon }+\lambda \right) \left( I+\sum\limits_{j=1}^{\infty }\Phi
_{\varepsilon j\lambda }\right) ^{-1}f=  \tag{2.20}
\end{equation}%
\ 

\begin{equation*}
\sum\limits_{j=1}^{\infty }\varphi _{j}\left( x\right) O_{\varepsilon
j\lambda }^{-1}{}\left[ I-K_{\varepsilon j\lambda }\right] ^{-1}\left(
I+\sum\limits_{j=1}^{\infty }\Phi _{\varepsilon j\lambda }\right) ^{-1}f.
\end{equation*}%
Then by using the above representation and by using Theorem A$_{1}$ we
obtain the estimate $\left( 2.3\right) $.

\textbf{Result 2.1. }Theorem 2.1 implies that the differential operator $%
O_{\varepsilon }$ has a resolvent $\left( O_{\varepsilon }+\lambda \right)
^{-1}$ for $\left\vert \arg \lambda \right\vert \leq \varphi ,$ and the
uniform estimate holds 
\begin{equation*}
\bigskip \sum\limits_{i=0}^{2}\left\vert \lambda \right\vert ^{1-\frac{i}{2}%
}\varepsilon ^{\frac{i}{2}}\left\Vert \frac{d^{i}}{dx^{i}}\left(
O_{\varepsilon }\mathbf{+}\lambda \right) ^{-1}\right\Vert _{L\left(
X\right) }+\left\Vert A\left( O_{\varepsilon }+\lambda \right)
^{-1}\right\Vert _{L\left( X\right) }\leq C.
\end{equation*}

\begin{center}
\textbf{3}. $\mathbf{R}$-\textbf{positive properties of the abstract
differential operator }
\end{center}

Result 2.1 implies that the operator $O$ is positive in $L_{p}\left( \sigma
;E\right) .$ In the following theorem we prove that this operator is $R$%
-positive of the operator $O$ in $L_{p}\left( \sigma ;E\right) .$

\textbf{Theorem 3.1. }Let all condition of Theorem 2.1 be satisfied. Then
the operator $O$ is $R$-positive in $L_{p}\left( \sigma ;E\right) .$

\textbf{Proof.} Consider first of all the problem with constant coefficents%
\begin{equation}
\ \varepsilon au^{\left( 2\right) }\left( x\right) +Au\left( x\right)
+\lambda u\left( x\right) =f\left( x\right) ,\text{ }x\in \sigma ,  \tag{3.1}
\end{equation}%
\begin{equation}
\sum\limits_{i=0}^{\mu _{1}}\varepsilon ^{\nu _{i}}\alpha _{i}u^{\left(
i\right) }\left( 0\right) =0,\text{ }\sum\limits_{i=0}^{\mu _{2}}\varepsilon
^{\nu _{i}}\beta _{i}u^{\left( i\right) }\left( 1\right) =0\text{, } 
\tag{3.2}
\end{equation}%
where $a$ is a complex number, $A$ is a linear operator in a Banach space $%
E, $ $\lambda $ is a complex parameter, $\varepsilon $ is a positive
parameter, $\nu _{i}=\frac{i}{2}+\frac{1}{2p},$ $\ \alpha _{i},$ $\beta _{i}$
are complex numbers, $\mu _{1},$ $\mu _{2}\in \left\{ 0,1\right\} .$

Consider the operator $O_{0}$ in $X$ generated by problem $\left( 3.1\right)
-\left( 3.2\right) $ for $\lambda =0$, i.e. 
\begin{equation*}
D\left( O_{0}\right) =W_{p}^{2}\left( \sigma ;E\left( A\right)
,E,L_{1},L_{2}\right) ,\text{ }O_{0}u=\varepsilon au^{\left( 2\right) }+Au.
\end{equation*}%
Since $A$ is a positive operator in $E,$ then in view of $\left[ \text{9,
Lemma 2.6}\right] $ there exists semigroups $U_{\varepsilon j\lambda }\left(
x\right) =e^{\varepsilon ^{\frac{1}{2}}x\omega _{1}A_{\lambda }^{\frac{1}{2}%
}}$for $\func{Re}\omega _{1}<0$, $U_{\varepsilon j\lambda }\left( x\right)
=e^{-\varepsilon ^{\frac{1}{2}}\left( b-x\right) \omega _{2}A_{\lambda }^{%
\frac{1}{2}}}$for $\func{Re}\omega _{2}>0\ $that are holomorphic for $x>0$
and strongly continuous for $x\geq 0$. By using a technique similar to that
applied in $\left[ \text{26},\text{ Lemma 5. 3. 2/1}\right] ,$ we obtain
that for $f$ $\in D\left( \sigma ;E\left( A\right) \right) $ the solution of
the equation $\left( 3.1\right) $ is represented as

\begin{equation}
u\left( x\right) =\dsum\limits_{j=1}^{2}U_{j\lambda }\left( x\right)
g_{k}+\int\limits_{\sigma }U_{0\lambda }\left( x-y\right) f\left( y\right) dy%
\text{, }g_{k}\in E,  \tag{3.3}
\end{equation}%
where%
\begin{equation*}
U_{\varepsilon 0\lambda }\left( x-y\right) =\left\{ 
\begin{array}{c}
-A_{\lambda }^{-\frac{1}{2}}U_{\varepsilon 1\lambda }\left( x-y\right) \text{%
, }x\geq y \\ 
A_{\lambda }^{-\frac{1}{2}}U_{\varepsilon 2\lambda }\left( y-x\right) \text{%
, }x\leq y%
\end{array}%
\right. .
\end{equation*}%
By taking into account the boundary conditions $\left( 3.2\right) ,$ we
obtain the following equation with respect to $g_{1},g_{2}$%
\begin{equation*}
\dsum\limits_{k=1}^{2}L_{k}\left( U_{\varepsilon j\lambda }\right)
g_{k}=L_{k}\left( \Phi _{\lambda }\right) \text{, }j=1,2,\text{ }\Phi
_{\lambda }=\int\limits_{\sigma }U_{\varepsilon 0\lambda }\left( x-y\right)
f\left( y\right) dy.
\end{equation*}%
By solving the above system and substituting it into $\left( 3.3\right) $ we
obtain the representation of the solution for problem $\left( 3.1\right)
-\left( 3.2\right) $:%
\begin{equation}
u\left( x\right) =\left[ O_{0}+\lambda \right] ^{-1}f=\int\limits_{\sigma
}G_{\varepsilon }\left( \lambda ,x,y\right) f\left( y\right) dy\text{, } 
\tag{3.4}
\end{equation}%
\begin{equation*}
G_{\varepsilon }\left( \lambda ,x,y\right)
=\sum_{k=1}^{2}\sum_{j=1}^{2}A_{\lambda }^{-\frac{1}{2}}B_{kj}\left( \lambda
\right) U_{j\lambda }\left( x\right) \tilde{U}_{kj\lambda }\left( x-y\right)
+U_{0\lambda }\left( x-y\right) ,
\end{equation*}%
where $B_{kj}\left( \lambda \right) $ are are uniformly bounded operators in 
$E$ and 
\begin{equation*}
\tilde{U}_{kj\lambda }\left( x-y\right) =\left\{ 
\begin{array}{c}
b_{kj}U_{\varepsilon k\lambda }\left( x-y\right) \text{, }x\geq y \\ 
\beta _{kj}U_{k\varepsilon \lambda }\left( y-x\right) \text{, }x\leq y%
\end{array}%
,b_{kj}\text{, }\beta _{kj}\in \mathbb{C}\mathbf{.}\right.
\end{equation*}%
Let at first, to show that the set $\Phi =\left\{ G_{\varepsilon }\left(
\lambda ,x,y\right) ;\lambda \in S\left( \varphi \right) \right\} $ is
uniformly $R$-bounded. By using the generalized Minkowcki's, Young
inequalities and by using of the holomorphic semigroups estimates $\left[ 9%
\right] $ we have the uniform estimate 
\begin{equation*}
\left\Vert G_{\varepsilon }\left( \lambda ,x,y\right) f\right\Vert _{X}\leq
C\sum_{k=1}^{2}\sum_{j=1}^{2}\left\{ \left\Vert A_{\lambda }^{-\frac{1}{2}%
}\right\Vert \left\Vert B_{kj}\left( \lambda \right) \right\Vert \left\Vert 
\tilde{U}_{kj\lambda }\left( x\right) f\right\Vert _{X}\right. +
\end{equation*}

\begin{equation*}
\left. \left\Vert U_{0\lambda }\left( x\right) f\right\Vert _{X}\right\}
\leq C\left\Vert f\right\Vert _{X}.
\end{equation*}%
Due to $R$-positivity of $A$, uniform boundedness of operators $B_{kj}\left(
\lambda \right) $ and in view of the Kahane's contraction principle and from
the product properties of the collection of $R$-bounded operators $\left[ 
\text{9, Lemma 3.5, Proposition 3.4}\right] $ we get that the sets%
\begin{eqnarray*}
b_{kj}\left( \lambda ,x,y\right) &=&\left\{ B_{kj}\left( \lambda \right)
A_{\lambda }^{-\frac{1}{2}}U_{\varepsilon j\lambda }\left( x\right) \left[
U_{\varepsilon k\lambda }\left( 1-y\right) +U_{\nu \lambda }\left( y\right) %
\right] :\lambda \in S_{\varphi }\right\} ,\text{ } \\
b_{0}\left( \lambda ,x,y\right) &=&\left\{ U_{\varepsilon 0\lambda }\left(
x-y\right) :\lambda \in S_{\varphi }\right\}
\end{eqnarray*}%
are uniformly $R$-bounded. Then by using the Kahane's contraction principle,
product and additional properties of the collection of $R$-bounded operators
and in view of $R$-boundedness of the sets $b_{kj}$, $b_{0},$ for all $%
u_{1,}u_{2},...,u_{\mu }\in F$, $\lambda _{1},\lambda _{2},...,\lambda _{\mu
}\in S\left( \varphi \right) $, and independent symmetric $\left\{
-1,1\right\} $-valued random variables $r_{i}\left( y\right) $, $%
i=1,2,...,\mu $, $\mu \in N$ \ we have the uniform estimate 
\begin{equation*}
\int\limits_{\Omega }\left\Vert \sum\limits_{i=1}^{\mu }r_{i}\left( y\right)
G_{\varepsilon }\left( \lambda _{i},x,y\right) u_{i}\right\Vert _{X}d\tau
\leq C\left\{ \sum_{k,j=1}^{2}\int\limits_{\Omega }\left\Vert
\sum\limits_{i=1}^{\mu }r_{i}\left( y\right) b_{kj}\left( \lambda
_{i},x,y\right) u_{i}\right\Vert _{X}d\tau \right.
\end{equation*}%
\begin{equation*}
+\int\limits_{\Omega }\left\Vert \sum\limits_{i=1}^{\mu }r_{i}\left(
y\right) b_{0}\left( \lambda _{i},x,y\right) u_{i}\right\Vert _{X}d\tau \leq
Ce^{\beta \left\vert \lambda \right\vert ^{\frac{1}{2}}\left\vert
x-y\right\vert }\int\limits_{\Omega }\left\Vert \sum\limits_{i=1}^{\mu
}r_{i}\left( y\right) u_{i}\right\Vert _{X}d\tau \text{, }\beta <0.
\end{equation*}%
This implies that 
\begin{equation*}
R\left\{ G_{\varepsilon }\left( \lambda ,x,y\right) :\lambda \in S_{\varphi
}\right\} \leq Ce^{\beta \left\vert \lambda \right\vert ^{\frac{1}{2}%
}\left\vert x-y\right\vert }\text{, }\beta <0\text{, }x,y\in \left(
0,b\right) .
\end{equation*}

By applying the $R$-bondedness property of kernel operators (see e.g. the
Proposition 4.12 in $\left[ \text{9}\right] $) and due to density of $%
D\left( \sigma ;E\left( A\right) \right) $ in $X$ ( see e.g.$\left[ \text{%
14, \S\ 2.2}\right] $ ) we \ get that the operator $O_{0}$ is uniformly $R$%
-positive in $X.$ From the representation $\left( 3.4\right) $ of solution
of problem $\left( 3.1\right) -\left( 3.2\right) $ it is easy to see that
the operator $\left( O_{0}+\lambda \right) ^{-1}$ can be expressed as a
linear combination of operators $O_{j\lambda }^{-1}$ like $\left(
O_{0}+\lambda \right) ^{-1}.$ Then, in view the representation $\left(
3.4\right) $ and by virtue of Kahane's contraction principle, product and
additional properties of the collection of $R$-bounded operators we obtain
that the operator $O_{0}$ is $R$-positive in $L_{p}\left( \sigma ;E\right) $.

Now, consider the problem $\left( 2.1\right) -\left( 2.2\right) .$ By virtue
of $\left( 2.20\right) $ from Theorem 2.1 we obtain that for$\ f\in
L_{p}\left( \sigma ;E\right) $ the BVP $\left( 2.1\right) -\left( 2.2\right) 
$ have a unique solution expressing in the form 
\begin{equation}
u\left( x\right) =\left( O_{\varepsilon }+\lambda \right)
^{-1}f=\sum\limits_{j=1}^{\infty }\varphi _{j}O_{\varepsilon j\lambda
}^{-1}{}\left[ I-K_{\varepsilon j\lambda }{}\right] ^{-1}g_{j}\left(
I+\sum\limits_{j=1}^{\infty }\Phi _{\varepsilon j\lambda }\right) ^{-1}f, 
\tag{3.5}
\end{equation}%
where $O_{\varepsilon j\lambda }=O_{\varepsilon j}+\lambda $\ are local
operators generated by BVPs with constant coefficients\ of type $\left(
2.16\right) $ and $K_{\varepsilon j\lambda },$ $\Phi _{\varepsilon j\lambda
} $ are uniformly bounded operators defined in the proof of the Theorem 2.1$%
. $ By virtue of the first part of this theorem, the operators $O_{j}$ are $%
R $-positive in $L_{p}\left( G_{j};E\right) $. Then by using the
representation $\left( 3.5\right) $ and by virtue of Kahane's contraction
principle, product and additional properties of the collection of $R$%
-bounded operators ( see e.g. $\left[ \text{9}\right] $ Lemma 3.5,
Proposition 3.4 ) we obtain the assertion.

\begin{center}
\textbf{4. Abstract Cauchy problem for parabolic equation on exterior domain}
\end{center}

Consider the following mixed problem for parabolic DOE equation with
parameter

\begin{equation*}
\ \frac{\partial u}{\partial t}+\varepsilon a\left( x\right) \frac{\partial
^{2}u}{\partial x^{2}}+\left[ A\left( x\right) +d\right] u+\varepsilon ^{%
\frac{1}{2}}A_{1}\left( x\right) \frac{\partial u}{\partial x}+A_{0}\left(
x\right) u=f\left( t,x\right) ,
\end{equation*}%
\begin{equation}
\sum\limits_{i=0}^{\mu _{1}}\varepsilon ^{\nu _{i}}\alpha _{i}u_{x}^{\left(
i\right) }\left( t,0\right) =0,\text{ }\sum\limits_{i=0}^{\mu
_{2}}\varepsilon ^{\nu _{i}}\beta _{i}u_{x}^{\left( i\right) }\left(
t,b\right) =0\text{, }  \tag{4.1}
\end{equation}

\begin{equation*}
u\left( 0,x\right) =0,\text{ }x\in \sigma ,\text{ }t\in \left( 0,T\right) ,
\end{equation*}%
where $\sigma =\left( -\infty ,\infty \right) \setminus \left[ 0,1\right] ,$ 
$\alpha _{i},$ $\beta _{i}$ are complex numbers, $\varepsilon $ is a
positive parameter, $\nu _{i}=\frac{i}{2}+\frac{1}{2p}$, $d$ is a positive
number$,$ $\mu _{k}\in \left\{ 0,1\right\} ,$ $A\left( .\right) $ and $%
A_{j}\left( .\right) $ are linear operator functions in a Banach space $E$
for $x\in \sigma .$

For $\mathbf{p=}\left( p,p_{1}\right) $, $\Delta _{+}=\left( 0,T\right)
\times \sigma ,$ $L_{\mathbf{p}}\left( \Delta _{T};E\right) $ will be
denoted the space of all $E$-valued $\mathbf{p}$-summable functions with
mixed norm (see e.g. $\left[ 6\right] $), i.e., the space of all measurable
functions $f$ defined on $\Delta _{T}$ for which 
\begin{equation*}
\left\Vert f\right\Vert _{L_{\mathbf{p}}\left( \Delta _{+}\right) }=\left(
\int\limits_{0}^{T}\left( \int\limits_{\sigma }\left\Vert f\left( t,x\right)
\right\Vert ^{p}dx\right) ^{\frac{p_{1}}{p}}dt\right) ^{\frac{1}{p_{1}}%
}<\infty .
\end{equation*}

Analogously, $W_{\mathbf{p}}^{2}\left( \sigma _{T},E\left( A\right)
,E\right) $ denotes the Sobolev space with corresponding mixed norm (see $%
\left[ 6\right] $ for scalar case)$.$

\textbf{Theorem 4.1. }Let the conditions of Theorem 2.1 hold for $\varphi >%
\frac{\pi }{2}$. Then for all $f\in L_{\mathbf{p}}\left( \sigma
_{T};E\right) $ and sufficiently large $d>0$ problem $\left( 4.1\right) $
has a unique solution belonging to $W_{\mathbf{p}}^{1,2}\left( \sigma
_{T};E\left( A\right) ,E\right) $ and the following coercive estimate holds 
\begin{equation*}
\left\Vert \frac{\partial u}{\partial t}\right\Vert _{L_{\mathbf{p}}\left(
\sigma _{T};E\right) }+\left\Vert \varepsilon \frac{\partial ^{2}u}{\partial
x^{2}}\right\Vert _{L_{\mathbf{p}}\left( \sigma _{T};E\right) }+\left\Vert
Au\right\Vert _{L_{\mathbf{p}}\left( \sigma _{+};E\right) }\leq C\left\Vert
f\right\Vert _{L_{\mathbf{p}}\left( \sigma _{T};E\right) }.
\end{equation*}%
\textbf{Proof. }The problem $\left( 4.1\right) $ can be express as the
following Cauchy problem 
\begin{equation}
\frac{du}{dt}+\left( O_{\varepsilon }+d\right) u=f\left( t\right) ,\text{ }%
u\left( 0\right) =0,  \tag{4.2}
\end{equation}%
where $O_{\varepsilon }$ denote the operator generated by $\left( 2.1\right)
-\left( 2.2\right) .$The Theorem 3.1 implies that the operator $%
O_{\varepsilon }$ is $R$-positive in $X=L_{p}\left( \sigma ;E\right) .$ By
virtue of $\left[ \text{24, \S 1.14}\right] ,$ the operator $O_{\varepsilon
} $ is a generator of an analytic semigroup in $X.$ Then applying $\left[ 
\text{9, Theorem 4.4}\right] $ we obtain that for $f\in L_{p_{1}}\left(
0,T;X\right) $ and sufficiently large $d>0$ problem $\left( 4.2\right) $ has
a unique solution belonging to $W_{p_{1}}^{1}\left( 0,T;D\left( O\right)
,X\right) $ and the following estimate holds 
\begin{equation*}
\left\Vert \frac{du}{dt}\right\Vert _{L_{p_{1}}\left( 0,T;X\right)
}+\left\Vert \left( O_{\varepsilon }+d\right) u\right\Vert _{L_{p_{1}}\left(
0,T;X\right) }\leq C\left\Vert f\right\Vert _{L_{p_{1}}\left( 0,T;X\right) }.
\end{equation*}%
Since $L_{p_{1}}\left( 0,T;X\right) =L_{\mathbf{p}}\left( \sigma
_{T};E\right) ,$ by Theorem 2.1 we have 
\begin{equation*}
\left\Vert \left( O_{\varepsilon }+d\right) u\right\Vert _{X}=\left\Vert
u\right\Vert _{W^{2}\left( \sigma ;E\left( A\right) ,E\right) }.
\end{equation*}%
These relations and the above estimate prove the hypothesis to be true.

\begin{center}
\ \textbf{5. Elliptic DOE on the moving domain\ }
\end{center}

Consider the BVP on the exterier moving domain $\sigma \left( s\right) =%
\mathbb{R}/\left[ 0,b\left( s\right) \right] $:%
\begin{equation}
au^{\left( 2\right) }\left( x\right) +Au\left( x\right) +A_{1}u^{\left(
1\right) }\left( x\right) +A_{0}u\left( x\right) =f\left( x\right) ,\text{ }%
x\in \sigma ,  \tag{5.1}
\end{equation}%
\begin{equation*}
\sum\limits_{i=0}^{\mu _{1}}\alpha _{i}u^{\left( i\right) }\left( 0\right)
=0,\text{ }\sum\limits_{i=0}^{\mu _{2}}\beta _{i}u^{\left( i\right) }\left(
b\left( s\right) \right) =0\text{, }
\end{equation*}%
\ where $\alpha _{i}$, $\beta _{i}$ are complex numbers, $a$ is a complex
valued function; $A=A\left( x\right) $ and $A_{j}=A_{j}\left( x\right) $ are
linear operators in a Banach space $E,$ the end point $b\left( s\right) $
depend on the parameter $s$ and $b(s)$ is a positive continues function on
compact domain $\Delta \subset \mathbb{R},$ $\mu _{k}\in \left\{ 0,1\right\}
.$

Theorem 2.1 implies the following:

\textbf{Proposition 5.1.} Assume the Condition 2.1 hold for $b=b\left(
s\right) $. Then, problem $\left( 5.1\right) $ has a unique solution $u\in $%
\ $W_{p}^{2}\left( \left( 0,b\right) ;E\left( A\right) ,E\right) $ for $f\in
L_{p}\left( 0,b;E\right) $ and sufficiently $d>0$. Moreover, the following
coercive uniform estimate holds

\begin{equation}
\left\Vert u^{2}\right\Vert _{L_{p}\left( 0,b;E\right) }+\left\Vert
Au\right\Vert _{L_{p}\left( 0,b;E\right) }\leq C\left\Vert f\right\Vert
_{L_{p}\left( 0,b;E\right) }.  \tag{5.2}
\end{equation}

\textbf{Proof. }Under the substitution $\tau =xb^{-1}(s)$ the problem $%
\left( 5.1\right) $ reduced to the following BVP in fixed domain $\left(
0,1\right) $:%
\begin{equation*}
b^{-2}\left( s\right) \tilde{a}\left( \tau \right) \tilde{u}^{\left(
2\right) }+\tilde{A}\left( \tau \right) \tilde{u}+\dsum%
\limits_{i=0}^{1}b^{-i}\left( s\right) \tilde{A}_{i}\left( \tau \right) 
\tilde{u}^{\left( i\right) }\left( \tau \right) =\tilde{f}\left( \tau
\right) \text{, }\tau \in \left( 0,1\right) ,
\end{equation*}%
\begin{equation*}
\dsum\limits_{i=0}^{1}b^{-i}\left( s\right) \alpha _{i}\tilde{u}^{\left(
i\right) }\left( 0\right) =0,\text{ }\dsum\limits_{i=0}^{1}b^{-i}\left(
s\right) \beta _{i}\tilde{u}^{\left( i\right) }\left( 1\right) =0,
\end{equation*}

where 
\begin{equation*}
\tilde{u}\left( \tau \right) =u\left( \tau b^{-1}\right) \text{, }\tilde{a}%
_{k}\left( \tau \right) =a_{k}\left( \tau b^{-1}\right) ,\text{ }\tilde{A}%
\left( \tau \right) =A\left( \left( \tau b^{-1}\right) \right) ,\text{ }
\end{equation*}%
\begin{equation*}
\tilde{A}_{i}\left( \tau \right) =A_{i}\left( \tau b^{-1}\right) ,\text{ }%
\tilde{f}\left( \tau \right) =f\left( \left( \tau b^{-1}\right) \right) .
\end{equation*}%
Then, by virtue of Theorem 2.1 we obtain the required assertion.

\begin{center}
\bigskip

\textbf{6.} \textbf{Nonlinear abstract elliptic problem in exterior domain}
\end{center}

Consider the following nonlinear elliptic problem 
\begin{equation}
a\left( x\right) u^{\left( 2\right) }\left( x\right) +B\left( x,u,u^{\left(
1\right) }\right) u=F\left( x,u,u^{\left( 1\right) }\right) +f\left(
x\right) ,\text{ }x\in \sigma ,  \tag{6.1}
\end{equation}

\begin{equation}
L_{1}u=\sum\limits_{i=0}^{\mu _{1}}\alpha _{i}u^{\left( i\right) }\left(
0\right) =0,\text{ }L_{2}u=\text{ }\sum\limits_{i=0}^{\mu _{2}}\beta
_{i}u^{\left( i\right) }\left( b\right) =0\text{,}  \tag{6.2}
\end{equation}%
where $g$ is $E$-valued given function, $a$ is a complex valued function, $%
\alpha _{i},$ $\beta _{i}$ $\ $are complex numbers, $\mu _{k}\in \left\{
0,1\right\} $, $\sigma =\mathbb{R}\setminus \left[ 0,b\right] .$

\bigskip In this section we will prove the existence and uniqueness of
maximal regular solution for the nonlinear problem $\left( 6.1\right)
-\left( 6.2\right) $. Let%
\begin{equation*}
U=\left( u_{0},u_{1}\right) ,\text{ }X=L_{p}\left( \sigma ;E\right) ,\text{ }%
Y=W_{p}^{2}\left( \sigma ;E\left( A\right) ,E\right) ,
\end{equation*}%
\begin{equation*}
E_{i}=\left( E\left( A\right) ,E\right) _{\theta _{i},p},\text{ }\theta _{i}=%
\frac{i+\frac{1}{p}}{2},\text{ }X_{0}=\prod\limits_{i=0}^{1}E_{i},
\end{equation*}

\bigskip \textbf{Remark 6.1. }By using J. Lions-I. Petree result (see e.g $%
\left[ \text{24, \S\ 1.8.}\right] $) we obtain that the embedding $D^{i}Y$ $%
\in E_{i}$ is continuous and there is a constant $C_{1}$ such that for $w\in
Y,$ $W=\left\{ w_{i}\right\} ,$ $w_{i}=D^{i}w\left( \cdot \right) ,$ $i=0,1,$%
\begin{equation*}
\left\Vert u\right\Vert _{\infty ,X_{0}}=\prod\limits_{i=0}^{1}\left\Vert
D^{i}w\right\Vert _{C(\bar{\sigma},E_{j})}=\sup_{x\in \left[ 0,b\right]
}\prod\limits_{i=0}^{1}\left\Vert D^{i}w\left( x\right) \right\Vert
_{E_{j}}\leq C_{1}\left\Vert w\right\Vert _{Y}\text{.}
\end{equation*}

For $r>0$ denote by $O_{r}$ the closed ball in $X_{0}$\ of radios $r$, i.e.

\begin{equation*}
O_{r}=\left\{ u\in X_{0}\text{, }\left\Vert u\right\Vert _{X_{0}}\leq
r\right\} .
\end{equation*}

Consider the linear problem, 
\begin{equation}
\ Lu=a\left( x\right) w^{\left( 2\right) }\left( x\right) +\left( A\left(
x\right) +d\right) w\left( x\right) =g\left( x\right) ,\text{ }  \tag{6.3}
\end{equation}

\begin{equation*}
L_{k}w=0,\text{ }k=1,2,
\end{equation*}%
where $A\left( x\right) $ is a linear operator in a Banach space $E$ for $%
x\in \sigma $, $L_{k}$ are boundary conditions defined by $\left( 6.1\right) 
$ and $d>0.$

Assume $E$ is a UMD space and $A\left( x\right) $ is uniformly $R$-positive
in $E,$ $A\left( 0\right) A^{-1}\left( y_{0}\right) =A\left( a\right)
A^{-1}\left( y_{0}\right) $. \ By virtue Theorem 2.1 and Proposition 5.1,
problem $\left( 6.3\right) $\ has a unique solution $w\in Y$ for all $g\in X$
and for sufficiently large $d>0$. Moreover, the following coercive estimate
holds 
\begin{equation*}
\left\Vert w\right\Vert _{Y}\leq M\left\Vert g\right\Vert _{X},
\end{equation*}%
where the constant $C_{0}$\ do not depend on $f\in X$ and $b\in \left(
0\right. \left. b_{0}\right] .$

Let $\omega _{1}=\omega _{1}\left( x\right) $, $\omega _{2}=\omega
_{2}\left( x\right) $ be roots of equation $a\left( x\right) \omega ^{2}+1=0$%
.

\textbf{Condition 6.1. }Assume the following satisfied:

(1 $a\in C\left( \bar{\sigma}\right) $, $\func{Re}\omega _{k}\neq 0$ and $%
\frac{\lambda }{\omega _{k}}\in S\left( \varphi \right) $ for $\lambda \in
S\left( \varphi \right) $, $0\leq \varphi <\pi ,$ $k=1,2$.$\ $a.e. $x\in 
\mathbb{\sigma };$

(2)\ $E$ is an UMD\ space$,$ $p\in \left( 1,\infty \right) $;

(3) $F:\bar{\sigma}\times X_{0}\rightarrow E$ is a measurable function for
each $u_{i}\in E_{i},$ $i=0,1$ and $F\left( x,U\right) \in X.$ Moreover, for
each $r>0$ there exists the positive functions $h_{k}\left( x\right) $ such
that 
\begin{equation*}
\left\Vert F\left( x,U\right) \right\Vert _{E}\leq h_{1}\left( x\right)
\left\Vert U\right\Vert _{X_{0}},
\end{equation*}%
\begin{equation*}
\left\Vert F\left( x,U\right) -F\left( x,\bar{U}\right) \right\Vert _{E}\leq
h_{2}\left( x\right) \left\Vert U-\bar{U}\right\Vert _{X_{0}},
\end{equation*}%
where $h_{k}\in L_{p}\left( \sigma \right) $ with 
\begin{equation*}
\left\Vert h_{k}\right\Vert _{L_{p}\left( \sigma \right) }<M^{-1}\text{, }%
k=1,2;
\end{equation*}%
and $U=\left\{ u_{0},u_{1}\right\} $, $\bar{U}=\left\{ \bar{u}_{0},\bar{u}%
_{1}\right\} $, $u_{i},$ $\bar{u}_{i}\in E_{i}$ and $U,$ $\bar{U}\in O_{r}.$

(4) there exist $\Phi _{i}\in E_{i}$, such that the operator $B\left( x,\Phi
\right) $ for $\Phi =\left\{ \Phi _{i}\right\} $ is $R$-positive in $E$
uniformly with respect to $x\in \left[ 0,b\right] ;$ $B\left( x,\Phi \right)
B^{-1}\left( x^{0},\Phi \right) \in C\left( \bar{\sigma};L\left( E\right)
\right) $; $B\left( x,0\right) =A\left( x\right) $;

(5) $B\left( x,U\right) $ for $x\in \left( 0,a\right) $ is a uniform
positive operator in a Banach space $E,$ where domain definition $D\left(
B\left( x,U\right) \right) $ does not depend on $x,U$ and $B$: $\sigma
\times X_{0}\rightarrow L\left( E\left( A\right) ,E\right) $ is continuous.
Moreover, for each $r>0$ there is a positive constant $L\left( r\right) $
such that

$\left\Vert \left[ B\left( x,U\right) -B\left( x,\bar{U}\right) \right]
\upsilon \right\Vert _{E}\leq L\left( r\right) \left\Vert U-\bar{U}%
\right\Vert _{X_{0}}\left\Vert A\upsilon \right\Vert _{E}$ for $x\in \sigma $%
, $U,\bar{U}\in O_{r}$ and $\upsilon \in D\left( B\left( x,U\right) \right)
. $

\textbf{Theorem 6.1. }Assume the Condition 6.1 holds. Then, there exist a
radius $0<r\leq r_{0}$ and $\delta $ $>0$ such that for each $f\in
L_{p}(\sigma ;E)$ with $\left\Vert f\right\Vert _{_{L_{p}\left( \sigma
;E\right) }}\leq \delta $ there is a unique solution $u\in W_{p}^{2}\left(
(\sigma ;E\left( A\right) ,E\right) $ of the problem $\left( 6.1\right)
-\left( 6.2\right) $ with $\left\Vert u\right\Vert _{W_{p}^{2}\left( \sigma
;E\left( A\right) ,E\right) }\leq r.$

\textbf{Proof.}\ We want to solve the problem $\left( 6.1\right) -\left(
6.2\right) $ locally by means of maximal regularity of the linear problem $%
(6.3)$ via the contraction mapping theorem. For this purpose, let $w$ be a
solution of the linear problem $(6.3).$ Consider a ball 
\begin{equation*}
O_{r}=\left\{ \upsilon \in Y,\text{ }L_{k}\left( \upsilon -w\right) =0,\text{
}\left\Vert \upsilon -w\right\Vert _{Y}\leq r\right\} .
\end{equation*}

Let $w\in Y$ be a solution of the problem $\left( 6.3\right) $ and 
\begin{equation*}
W=\left\{ w,w^{\left( 1\right) }\right\} .
\end{equation*}%
Given $\upsilon \in B_{r}$ solve the linear problem%
\begin{equation}
a\left( x\right) u^{\left( 2\right) }\left( x\right) +B\left( x,0\right)
u\left( x\right) +du=F\left( x,\upsilon \right) +  \notag
\end{equation}%
\begin{equation}
\left[ B\left( x,0\right) -B\left( x,V\right) \right] \upsilon \left(
x\right) +f\left( x\right) ,\text{ }L_{k}u=0,\text{ }k=1,2,  \tag{6.4}
\end{equation}%
where 
\begin{equation*}
V=\left( \upsilon ,\upsilon ^{\left( 1\right) }\right) \text{, }\upsilon \in
Y.
\end{equation*}%
Consider the function 
\begin{equation}
g\left( x\right) =\left[ B\left( x,0\right) -B\left( x,V\right) \right]
\upsilon \left( x\right) +\text{ }F\left( x,V\right) +f\left( x\right) . 
\tag{6.5}
\end{equation}

\bigskip Let first of all, we show that $g\in X$ and $\left\Vert
g\right\Vert _{X}\leq M^{-1}r$ for $\upsilon \in Y,$ $\left\Vert \upsilon
\right\Vert _{Y}\leq r.$ Indeed, by Remark 6.1 $V\in C\left( \bar{\sigma}%
;X_{0}\right) $, one has 
\begin{equation*}
B\left( x,0\right) -B\left( x,V\right) \in C\left( \bar{\sigma};L\left(
E\left( A\right) ,E\right) \right) .
\end{equation*}%
Hence, by assumption (3), $g$ is measurable and 
\begin{equation*}
\left\Vert g\left( x\right) \right\Vert _{E}\leq L\left( r\right) \left\Vert
V\right\Vert _{X_{0}}+h_{1}\left( x\right) \left\Vert V\right\Vert
_{X_{0}}+\left\Vert f\left( x\right) \right\Vert _{E}
\end{equation*}%
for a.e. $x\in \sigma $. Then, by using the Remark 6.1 and by chousing $%
\delta $ we obtain%
\begin{equation*}
\left\Vert g\right\Vert _{X}\leq rL\left( r\right) \left\Vert \upsilon
\right\Vert _{X}+r\left\Vert h_{1}\right\Vert _{L_{p}}+\left\Vert
f\right\Vert _{X}\leq r^{2}L\left( r\right) +r\left\Vert h_{1}\right\Vert
_{L_{p}}+\delta \leq M^{-1}r.
\end{equation*}%
Define a map $Q$ on $O_{r}$ by%
\begin{equation*}
Q\upsilon =w,
\end{equation*}%
where $w$ is a solution of the problem $\left( 6.3\right) $ with $g$ defined
by $\left( 6.4\right) .$ We want to show that $Q\left( B_{r}\right) \subset
B_{r}$ and that $Q$ is a contraction operator in $Y$ provided $\delta $ is
sufficiently small, and $r$ is chosen properly. For this aim, by using
maximal regularity properties of the problem $\left( 6.3\right) $ we have 
\begin{equation*}
\left\Vert Q\upsilon -w\right\Vert _{Y}=\left\Vert u-w\right\Vert _{Y}\leq
M\left\{ \left\Vert F\left( x,V\right) -F\left( x,0\right) \right\Vert
_{X}+\right.
\end{equation*}%
\begin{equation*}
\left. \left\Vert \left[ B\left( x,0\right) -B\left( x,V\right) \right]
\upsilon \right\Vert _{X}\right\} .
\end{equation*}

By assumption (3) for $\upsilon \in O_{r}$\ we get 
\begin{equation*}
\left\Vert F\left( x,V\right) -F\left( x,0\right) \right\Vert _{X}\leq
\left\Vert h_{2}\right\Vert _{L_{p}\left( \sigma \right) }\left\Vert
V\right\Vert _{X_{0}}.
\end{equation*}

By assumptions (4), (5) and Remark 6.2, for $\upsilon \in O_{r}$ and $%
W=\left( w,w^{\left( 1\right) }\right) $, $w\in Y$\ we have 
\begin{equation*}
\left\Vert \left[ B\left( x,0\right) \upsilon -B\left( x,V\right) \right]
\upsilon \right\Vert _{X}\leq \sup\limits_{x\in \bar{\sigma}}\left\{
\left\Vert \left[ B\left( x,0\right) -B\left( x,W\right) \right] \upsilon
\right\Vert _{L\left( X_{0},X\right) }\right.
\end{equation*}%
\begin{equation*}
+\left. \left\Vert \left[ B\left( x,W\right) -B\left( x,V\right) \right]
\upsilon \right\Vert _{L\left( X_{0},X\right) }\left\Vert \upsilon
\right\Vert _{Y}\right\} \leq
\end{equation*}%
\begin{equation*}
L\left( r\right) \left[ \left\Vert W\right\Vert _{X_{0}}\left\Vert A\upsilon
\right\Vert _{X}+\left\Vert \upsilon -w\right\Vert _{\infty ,X_{0}}\right] %
\left[ \left\Vert \upsilon -w\right\Vert _{Y}+\left\Vert w\right\Vert _{Y}%
\right] \leq
\end{equation*}%
\begin{equation*}
rL\left( r\right) \left\{ \left[ \left\Vert W\right\Vert _{X_{0}}\left\Vert
\upsilon \right\Vert _{Y}+C_{1}\left\Vert \upsilon -w\right\Vert _{Y}\right]
\left. +L\left( r\right) \left\Vert w\right\Vert _{Y}\right\} \right. .
\end{equation*}

By chousing $r$ and $b\in \left( 0\right. \left. b_{0}\right] $ so that $%
\left\Vert w\right\Vert _{Y}<\delta _{a}$ by assumptions (3)-(5) we obtain
from the above inequalities 
\begin{equation*}
\left\Vert Q\upsilon -w\right\Vert _{Y}\leq r+r^{2}L\left( r\right)
\left\Vert W\right\Vert _{X_{0}}+r^{2}L\left( r\right) C_{1}+rL\left(
r\right) \left\Vert w\right\Vert _{Y}<r.
\end{equation*}

That is the operator $Q$ maps $B_{r}$ into itself, i.e. 
\begin{equation*}
Q\left( B_{r}\right) \subset B_{r}.
\end{equation*}

Let $u_{1}=Q\left( \upsilon _{1}\right) $ and $u_{2}=Q\left( \upsilon
_{2}\right) .$ Then $u_{1}-u_{2}$ is a solution of the problem%
\begin{equation}
a\left( x\right) u^{\left( 2\right) }\left( x\right) +A\left( x\right)
u\left( x\right) +du=F\left( x,\upsilon _{1}\right) -  \notag
\end{equation}%
\begin{equation*}
F\left( x,\upsilon _{1}\right) +\left[ B\left( x,\upsilon _{2}\right)
-B\left( x,0\right) \right] \left[ \upsilon _{1}\left( x\right) -\upsilon
_{2}\left( x\right) \right] -
\end{equation*}%
\begin{equation*}
\left[ B\left( x,\upsilon _{1}\right) -B\left( x,\upsilon _{2}\right) \right]
\upsilon _{1}\left( x\right) ,\text{ }L_{k}u=0,\text{ }k=1,2.
\end{equation*}%
\ In a similar way, by using the assumption (5) we obtain%
\begin{equation*}
\left\Vert u_{1}-u_{2}\right\Vert _{Y}\leq C_{0}\left\{ rL\left( r\right)
\left\Vert \upsilon _{1}-\upsilon _{2}\right\Vert _{X}\right. +L\left(
r\right) \left\Vert \upsilon _{1}-\upsilon _{2}\right\Vert _{Y}\left\Vert
\upsilon _{1}\right\Vert _{X}
\end{equation*}%
\begin{equation*}
+\left. \left\Vert h_{2}\right\Vert _{L_{p}}\left\Vert \upsilon
_{1}-\upsilon _{2}\right\Vert _{Y}\right\} \leq C_{0}\left[ 2rL\left(
r\right) +\left\Vert h_{2}\right\Vert _{L_{p}}\right] \left\Vert \upsilon
_{1}-\upsilon _{2}\right\Vert _{Y}.
\end{equation*}

Thus $Q$ is a strict contraction. Eventually, the contraction mapping
principle implies a unique fixed point of $Q$ in $O_{r}$ which is the unique
strong solution 
\begin{equation*}
u\in Y=W_{p}^{2}\left( \sigma ;E\left( A\right) ,E\right) .
\end{equation*}

\begin{center}
\textbf{7. Exterior BVP for elliptic equations }

\ \ \ \ \ \ \ \ \ \ \ \ \ \ \ \ \ \ \ \ \ \ \ \ \ \ \ \ \ \ \ \ \ \ \ \ \ \
\ \ \ \ \ \ \ \ \ \ 
\end{center}

The regularity property of BVP for elliptic equations\ were studied e.g. in $%
\left[ 1\right] ,$ $\left[ 9\right] ,$ $\left[ 26\right] $. Let $\Omega
=\sigma \times G$, where $\sigma =\mathbb{R}\setminus \left[ 0,b\right] ,$ $%
G\subset R^{n},$ $n\geq 2$ is a bounded domain with $\left( n-1\right) $%
-dimensional boundary $\partial G$. Let us consider the following BVP for
elliptic equation with parameter 
\begin{equation}
\ Lu=\varepsilon a\left( x\right) D_{x}^{2}u\left( x,y\right)
+\sum\limits_{\left\vert \alpha \right\vert \leq 2m}b_{\alpha }\left(
x\right) a_{\alpha }\left( y\right) D_{y}^{\alpha }u\left( x,y\right) + 
\notag
\end{equation}%
\begin{equation}
+\sum\limits_{i=0}^{1}\dsum\limits_{\left\vert \beta \right\vert \leq \mu
_{i}}a_{_{i\beta }}\left( x,y\right) D_{x}^{i}D_{y}^{\beta }u\left(
x,y\right) +du\left( x,y\right) =f,\text{ }x\in \sigma \text{, }y\in G, 
\tag{7.1}
\end{equation}%
\begin{equation}
\sum\limits_{i=0}^{\eta _{1}}\varepsilon ^{\nu _{i}}\alpha _{i}u^{\left(
i\right) }\left( 0,y\right) =0,\text{ }\sum\limits_{i=0}^{\eta
_{2}}\varepsilon ^{\nu _{i}}\beta _{i}u^{\left( i\right) }\left( 0,y\right)
=0,\text{ for a.e. }y\in G,  \tag{7.2}
\end{equation}

\begin{equation}
B_{j}u=\sum\limits_{\left\vert \beta \right\vert \leq m_{j}}\ b_{j\beta
}\left( y\right) D_{y}^{\beta }u\left( x,y\right) =0\text{, }x\in \sigma ,%
\text{ }y\in \partial \Omega ,\text{ }j=1,2,...,m,  \tag{7.3}
\end{equation}%
where $\eta _{k}\in \left\{ 0,1\right\} ,$ $\alpha _{i},$ $\beta _{i}$ are
complex numbers, $\ \varepsilon $ is a positive parameter$,\nu _{i}=\frac{i}{%
2}+\frac{1}{2p}$, $d>0$,%
\begin{equation*}
D_{x}^{k}=\frac{\partial ^{k}}{\partial x^{k}},\text{ }D_{j}=-i\frac{%
\partial }{\partial y_{j}},\text{ }D_{y}=\left( D_{1,}...,D_{n}\right) ,%
\text{ }y=\left( y_{1},...,y_{n}\right)
\end{equation*}%
$\ $and $a,$ $a_{\alpha },$ $b_{\alpha },$ $a_{i\beta },$ $b_{j\beta }$ are
the complex valued functions, $\mu _{i}<2m.$ Let $\mathbf{p=}\left(
p_{1},p\right) .$

Let $\xi ^{\prime }=\left( \xi _{1},\xi _{2},...,\xi _{n-1}\right) \in
R^{n-1},$ $\alpha ^{\prime }=\left( \alpha _{1},\alpha _{2},...,\alpha
_{n-1}\right) \in Z^{n}$ and 
\begin{equation*}
\text{ }A\left( y_{0},\xi ^{\prime },D_{y}\right) =\sum\limits_{\left\vert
\alpha ^{\prime }\right\vert +j\leq 2m}a_{\alpha ^{\prime }}\left(
y_{0}\right) \xi _{1}^{\alpha _{1}}\xi _{2}^{\alpha _{2}}...\xi
_{n-1}^{\alpha _{n-1}}D_{y}^{j}\text{ for }y_{0}\in \bar{G}
\end{equation*}%
\begin{equation*}
B_{j}\left( y_{0},\xi ^{\prime },D_{y}\right) =\sum\limits_{\left\vert \beta
^{\prime }\right\vert +j\leq m_{j}}b_{j\beta ^{\prime }}\left( y_{0}\right)
\xi _{1}^{\beta _{1}}\xi _{2}^{\beta _{2}}...\xi _{n-1}^{\beta
_{n-1}}D_{y}^{j}\text{ for }y_{0}\in \partial G.
\end{equation*}

Let $Q$ denote the differential operator in $L_{\mathbf{p}}\left( \Omega
\right) $ generated by BVP $\left( 7.1\right) -\left( 7.3\right) .$

\textbf{Theorem 5.1}. Let the following conditions be satisfied:

\bigskip (1) $a\in C\left( \bar{\sigma}\right) $, $\func{Re}\omega _{k}\neq
0 $ and $\frac{\lambda }{\omega _{k}}\in S\left( \varphi \right) $ for $%
\lambda \in S\left( \varphi \right) $, $0\leq \varphi <\pi ,$ $k=1,2$.$\ $%
a.e. $x\in \mathbb{\sigma },\ b_{\alpha }\in C\left( \sigma \right) $, $%
a_{\alpha }\in C\left( \bar{\Omega}\right) $ for each $\left\vert \alpha
\right\vert =2m$ and $a_{\alpha }\in L_{\infty }\left( \Omega \right) $ for
each $\left\vert \alpha \right\vert <2m$;

(2) $b_{j\beta }\in C^{2m-m_{j}}\left( \partial \Omega \right) $ for each $j$%
, $\beta $ and $\ m_{j}<2m$, $\sum\limits_{j=1}^{m}b_{j\beta }\left(
y^{^{\shortmid }}\right) \sigma _{j}\neq 0,$ for $\left\vert \beta
\right\vert =m_{j},$ $y^{^{\shortmid }}\in \partial G,$ where $\sigma
=\left( \sigma _{1},\sigma _{2},...,\sigma _{n}\right) \in R^{n}$ is a
normal to $\partial G$ $;$

(3) for $y\in \bar{\Omega}$, $\xi \in R^{n}$, $\lambda \in S\left( \varphi
_{0}\right) $, $\left\vert \xi \right\vert +\left\vert \lambda \right\vert
\neq 0$ let $\lambda +$ $\sum\limits_{\left\vert \alpha \right\vert
=2m}a_{\alpha }\left( y\right) \xi ^{\alpha }\neq 0$;

(4) for each $y_{0}\in \partial \Omega $ local BVP in local coordinates
corresponding to $y_{0}$%
\begin{equation*}
\lambda +A\left( y_{0},\xi ^{\prime },D_{y}\right) \vartheta \left( y\right)
=0,
\end{equation*}

\begin{equation*}
B_{j}\left( y_{0},\xi ^{\prime },D_{y}\right) \vartheta \left( 0\right)
=h_{j}\text{, }j=1,2,...,m
\end{equation*}%
has a unique solution $\vartheta \in C_{0}\left( \mathbb{R}_{+}\right) $ for
all $h=\left( h_{1},h_{2},...,h_{n}\right) \in \mathbb{C}^{n}$ and for $\xi
^{\prime }\in R^{n-1}.$

Then;

(a) problem $\left( 7.1\right) -\left( 7.3\right) $ has a unique solution $%
u\in W_{\mathbf{p}}^{2,2m}\left( \Omega \right) $ for\ $f\in L_{\mathbf{p}%
}\left( \Omega \right) $ and sufficiently large $d>0.$ Moreover, the uniform
coercive estimate holds 
\begin{equation*}
\varepsilon \left\Vert D_{x}^{2}u\right\Vert _{L_{\mathbf{p}}\left( \Omega
\right) }+\sum\limits_{\left\vert \alpha \right\vert \leq 2m}\left\Vert
D^{\alpha }u\right\Vert _{L_{\mathbf{p}}\left( \Omega \right) }\leq
C\left\Vert f\right\Vert _{L_{\mathbf{p}}\left( \Omega \right) };
\end{equation*}

(b) the operator $Q$ is $R$-positive in $L_{\mathbf{p}}\left( \Omega \right)
.$

\ \textbf{Proof. }Let us consider operators $A$ and $A_{i}\left( x\right) $
in $E=L_{p_{1}}\left( G\right) $ that are defined by the equalities 
\begin{equation*}
D\left( A\right) =\left\{ u\in W_{p_{1}}^{2m}\left( G\right) \text{, }%
B_{j}u=0,\text{ }j=1,2,...,m\text{ }\right\} ,\ Au=\sum\limits_{\left\vert
\alpha \right\vert \leq 2m}a_{\alpha }\left( y\right) D_{y}^{\alpha }u\left(
y\right) ,
\end{equation*}

\begin{equation*}
\text{ }A_{i}u=\sum\limits_{\left\vert \beta \right\vert \leq \mu
_{i}}a_{i\beta }\left( x,y\right) D_{y}^{\beta }u\left( y\right) \text{, }%
i=0,1.
\end{equation*}%
Then the problem $\left( 7.1\right) -\left( 7.3\right) $ can be rewritten as
the problem $\left( 2.1\right) -\left( 2.2\right) $, where $u\left( x\right)
=u\left( x,.\right) ,$ $f\left( x\right) =f\left( x,.\right) $,\ $x\in
\sigma $ are the functions with values in\ $E=L_{p_{1}}\left( G\right) $. By
virtue of $\left[ \text{2, Theorem 4.5.2}\right] $ ) the space $%
E=L_{p_{1}}\left( G\right) ,$ $p_{1}\in \left( 1,\infty \right) $ satisfies
the multiplier condition. By virtue of $\left[ \text{9, Theorem 8.2}\right] $
operator $A+\mu $ for sufficiently large $\mu >0$ is $R$-positive in $%
L_{p_{1}}$. Moreover, (1) and (2) implies the (3) condition of Theorem 2.1,
i.e., conditions (1)- (3) of Theorem 2.1 are fulfilled. It is known that the
embedding $W_{p_{1}}^{2m}\left( G\right) \subset L_{p_{1}}\left( G\right) $
\ is compact (see e.g. $\left[ \text{24, \S\ 3}\right] ,$Theorem 3. 2. 5 ).
Using interpolation properties of Sobolev spaces $\left[ \text{24, \S\ 4}%
\right] $ we obtain that the condition (4) of Theorem 2.1 is satisfied.
Hence, all hypotheses of Theorem 2.1 are valid and the assertion of (a)
holds. Then the Theorem 3.1 implies the assertion (b).

\begin{center}
\textbf{8. The system of parabolic equation of arbitrary number on exterior
domain }
\end{center}

Consider the Cauchy problem for the system of parabolic equation of
arbitrary number%
\begin{equation*}
\ \frac{\partial u_{j}}{\partial t}+\varepsilon a\left( x\right) \frac{%
\partial ^{2}u_{j}}{\partial x^{2}}+\left[ a_{j}\left( x\right) +d\right]
u+\dsum\limits_{i=0}^{1}b_{ij}\left( x\right) \frac{\partial u_{j}}{\partial
x}+du=f\left( t,x\right) ,
\end{equation*}%
\begin{equation}
\sum\limits_{i=0}^{\mu _{1}}\varepsilon ^{\nu _{i}}\alpha _{i}\frac{\partial
u_{j}}{\partial x}\left( t,0\right) =0,\text{ }\sum\limits_{i=0}^{\mu
_{2}}\varepsilon ^{\nu _{i}}\beta _{i}\frac{\partial u_{j}}{\partial x}%
\left( t,1\right) =0\text{, }  \tag{8.1}
\end{equation}

\begin{equation*}
u\left( 0,x\right) =0,\text{ }x\in \sigma ,\text{ }t\in \left( 0,T\right) ,%
\text{ }j=1,2,...,\text{ }N\in \mathbb{N},
\end{equation*}%
where $a\left( .\right) ,$ $a_{j}\left( .\right) ,$ $b_{ij}\left( x\right) $
are complex valued functions, $\alpha _{i},$ $\beta _{i}$ are complex
numbers, $\varepsilon $ is a small positive parameter, $\nu _{i}=\frac{i}{2}+%
\frac{1}{2p},\ d$ is a positive number$,$ $\mu _{k}\in \left\{ 0,1\right\} , 
$ $\sigma =\mathbb{R}\setminus \left[ 0,1\right] .$

Let $\mathbf{p=}\left( p,p_{1}\right) $, $\Delta _{+}=\left( 0,T\right)
\times \sigma $ and $L_{\mathbf{p}}\left( \Delta _{T}\right) =L_{\mathbf{p}%
}\left( \Delta _{T};\mathbb{C}\right) $ will be denoted the space of all
complex-valued functions with mixed norm i.e., the space of all measurable
functions $f$ defined on $\Delta _{T}$ for which 
\begin{equation*}
\left\Vert f\right\Vert _{L_{\mathbf{p,\gamma }}\left( \Delta _{+}\right)
}=\left( \int\limits_{0}^{T}\left( \int\limits_{\sigma }\left\vert f\left(
t,x\right) \right\vert ^{p}dx\right) ^{\frac{p_{1}}{p}}dt\right) ^{\frac{1}{%
p_{1}}}<\infty .
\end{equation*}%
Analogously, $W_{\mathbf{p}}^{2}\left( \Delta _{T}\right) $ denotes the
Sobolev space with corresponding mixed norm (see e.g. $\left[ 6\right] $)$.$

\bigskip\ Let $E=l_{q}$ and $A\left( x\right) =\left[ \delta
_{ij}a_{i}\left( x\right) \right] ,$ $A_{i}\left( x\right) =\left[
b_{ij}\left( x\right) \right] $ are diagonal matrices in $l_{q},$ where $i$, 
$j=1,2,...N$, $\delta _{ij}=1$ for $i=j$ and $\delta _{ij}=0$ and 
\begin{equation*}
\text{ }l_{q}\left( A\right) =\left\{ u\in l_{q},\left\Vert u\right\Vert
_{l_{q}\left( A\right) }=\left\Vert Au\right\Vert _{l_{q}}=\right.
\end{equation*}

\begin{equation*}
\left( \sum\limits_{j=1}^{N}\left\vert \left( Au\right) _{j}\right\vert
^{q}\right) ^{\frac{1}{q}}=\left. \left( \left\vert
\sum\limits_{j=1}^{N}a_{j}u_{j}\right\vert ^{q}\right) ^{\frac{1}{q}}<\infty
\right\} ,
\end{equation*}%
\begin{equation*}
\text{ }u=\left\{ u_{j}\right\} ,\text{ }Au=\left\{
\dsum\limits_{j=1}^{N}a_{j}u_{j}\right\} ,\text{ }j=1,2,...N.
\end{equation*}

\textbf{Condition 8.1.} Assume the following conditions are satisfied;

(1)\textbf{\ } $a\in C\left( \bar{\sigma}\right) $, $\func{Re}\omega
_{k}\neq 0$ and $\frac{\lambda }{\omega _{k}}\in S\left( \varphi \right) $
for $\lambda \in S\left( \varphi \right) $, $k=1,2$, $a_{j}\in C\left( \bar{%
\sigma}\right) $ and $a_{j}\left( x\right) \in S\left( \varphi \right) ,$ $%
x\in \sigma ,$ $0\leq \varphi <\pi ;$

(2) $b_{ij}\in L_{\infty }\left( 0,1\right) ,$ $\left\vert b_{ij}\left(
x\right) \right\vert \leq C\left\vert a_{j}^{1-\frac{i}{2}-\delta
_{i}}\left( x\right) \right\vert $ for $0<\delta _{i}<1-\frac{i}{2}$ and
a.e. $x\in \sigma ;$

(5) $p,$ $q\in \left( 1,\infty \right) $ and $\dprod\limits_{j=1}^{N}\left%
\vert a_{j}\left( x\right) \right\vert <\infty $ for a.e. $x\in \sigma .$

Let%
\begin{equation*}
f\left( x\right) =\left\{ f_{j}\left( x\right) \right\} _{1}^{N}\text{, }%
u=\left\{ u_{j}\left( x\right) \right\} _{1}^{N}.
\end{equation*}

\textbf{Theorem 8.1. }Assume\textbf{\ }Condition 8.1 are satisfied.\textbf{\ 
}Then for $f\left( x\right) \in L_{p}\left( \Delta _{+};l_{q}\right) $ and
for sufficiently large $d$ problem $\left( 8.1\right) $ has a unique
solution $u$ that belongs to the space $W_{p}^{1,2}\left( \Delta
_{+};l_{q}\left( A\right) ,l_{q}\right) $ and the following coercive
estimate holds 
\begin{equation*}
\left[ \int\limits_{\Delta _{+}}\left( \sum\limits_{j=1}^{N}\left\vert \frac{%
\partial u_{j}}{\partial t}\right\vert ^{q}\right) ^{\frac{p}{q}}dx\right] ^{%
\frac{1}{p}}+\left[ \int\limits_{\Delta _{+}}\left(
\sum\limits_{j=1}^{N}\left\vert \frac{\partial ^{2}u_{j}}{\partial x^{2}}%
\right\vert ^{q}\right) ^{\frac{p}{q}}dx\right] ^{\frac{1}{p}}
\end{equation*}

\begin{equation*}
\left[ \int\limits_{\Delta _{+}}\left( \left\vert
\sum\limits_{j=1}^{N}a_{j}u_{j}\right\vert ^{q}\right) ^{\frac{p}{q}}dx%
\right] ^{\frac{1}{p}}\leq C\left[ \int\limits_{\Delta _{+}}\left(
\sum\limits_{i=1}^{N}\left\vert f_{i}\left( x\right) \right\vert ^{q}\right)
^{\frac{p}{q}}dx\right] ^{\frac{1}{p}}.
\end{equation*}

\textbf{Proof.} Let first all of, we suppose $N<\infty .$ Then det $A\left(
x\right) =\dprod\limits_{j=1}^{N}a_{j}\left( x\right) .$

It is easy to see that 
\begin{equation*}
B\left( \lambda \right) =\lambda \left( A+\lambda \right) ^{-1}=\frac{%
\lambda }{D\left( \lambda \right) }\left[ A_{ji}\left( \lambda \right) %
\right] \text{, }i\text{, }j=1,2,...N,
\end{equation*}%
where $D\left( \lambda \right) =\dprod\limits_{j=1}^{N}\left( a_{j}\left(
x\right) +\lambda \right) ^{-1}$, $A_{ji}\left( \lambda \right) $ are
entries of the corresponding adjoint matrix of $A+\lambda I.$ By using the
(1) assumption it is clear to see that the matrix $A$ generates a positive
operator in $l_{q}.$ For all $u_{1,}u_{2},...,u_{\mu }\in l_{q}$, $\lambda
_{1},\lambda _{2},...,\lambda _{\mu }\in \mathbb{C}$ and independent
symmetric $\left\{ -1,1\right\} $-valued random variables $r_{k}\left(
y\right) $, $k=1,2,...,\mu ,$ $\mu \in \mathbb{N}$\ we have 
\begin{equation*}
\int\limits_{0}^{1}\left\Vert \sum\limits_{k=1}^{\mu }r_{k}\left( y\right)
B\left( \lambda _{k}\right) u_{k}\right\Vert _{l_{q}}^{q}dy\leq
\end{equation*}%
\begin{equation*}
C\left\{ \int\limits_{0}^{1}\sum\limits_{j=1}^{N}\left\vert
\sum\limits_{k=1}^{\mu }\sum\limits_{j=1}^{N}\frac{\lambda _{k}}{D\left(
\lambda _{k}\right) }A_{ji}\left( \lambda _{k}\right) r_{k}\left( y\right)
u_{ki}\right\vert ^{q}dy\right. \leq
\end{equation*}%
\begin{equation}
\sup\limits_{k,i}\sum\limits_{j=1}^{N}\left\vert \frac{\lambda _{k}}{D\left(
\lambda _{k}\right) }A_{ji}\left( \lambda _{k}\right) \right\vert
^{q}\int\limits_{\Omega }\left\vert \sum\limits_{k=1}^{\mu }r_{k}\left(
y\right) u_{kj}\right\vert ^{q}dy.  \tag{8.2}
\end{equation}%
Since $A$ is symmetric and positive definite, we have%
\begin{equation}
\sup\limits_{k,i}\sum\limits_{j=1}^{N}\left\vert \frac{\lambda _{k}}{D\left(
\lambda _{k}\right) }A_{ji}\left( \lambda _{k}\right) \right\vert ^{q}\leq C.
\tag{8.3}
\end{equation}%
From $\left( 8.2\right) $ and $\left( 8.3\right) $ we get 
\begin{equation*}
\int\limits_{0}^{1}\left\Vert \sum\limits_{k=1}^{\mu }r_{k}\left( y\right)
B\left( \lambda _{k}\right) u_{k}\right\Vert _{l_{q}}^{q}dy\leq
C\int\limits_{0}^{1}\left\Vert \sum\limits_{k=1}^{\mu }r_{k}\left( y\right)
u_{k}\right\Vert _{l_{q}}^{q}dy,
\end{equation*}%
i.e., the operator $A$ is $R$-positive in $l_{q}.$

Let $N=\infty $, then we define determinant of infinite dimensional matrix $%
A $ as:%
\begin{equation*}
\text{det }A=\lim\limits_{n\rightarrow \infty
}\dprod\limits_{j=1}^{n}a_{j}<\infty .
\end{equation*}

The resolvent set $R\left( A\right) $ of the infinite dimensional matrix $A$
is defined as:%
\begin{equation*}
R\left( A\right) =\left\{ \lambda \in \mathbb{C}\text{, }\lim\limits_{n%
\rightarrow \infty }\dprod\limits_{j=1}^{n}\left( a_{j}+\lambda \right)
^{-1}<\infty \right\} .
\end{equation*}%
In a similar way we obtain that 
\begin{equation*}
B\left( \lambda \right) =\lambda \left( A+\lambda \right) ^{-1}=\frac{%
\lambda }{D\left( \lambda \right) }\left[ A_{ji}\left( \lambda \right) %
\right] \text{, }i\text{, }j=1,2,...N,
\end{equation*}%
where $D\left( \lambda \right) =\lim\limits_{n\rightarrow \infty
}\dprod\limits_{j=1}^{n}\left( a_{j}+\lambda \right) ^{-1}$ and $%
A_{ji}\left( \lambda \right) $ are entries of the corresponding adjoint
matrix of $A+\lambda .$ By reasoning as the above and by taking limit when $%
n\rightarrow \infty $ we obtain that the matrix $A$ generates $R-$positive
operator in $l_{q}$ also for $N=\infty .$ From the Theorem 3.1 we obtain
that problem $\left( 8.1\right) $ has a unique solution $u\in $ $W_{\mathbf{p%
}}^{1,2}\left( \Delta _{+};l_{q}\left( A\right) ,l_{q}\right) $ for $f\in L_{%
\mathbf{p}}\left( \Delta _{+};l_{q}\right) $ and the following uniform
estimate holds%
\begin{equation*}
\left\Vert \frac{\partial u}{\partial t}\right\Vert _{L_{\mathbf{p}}\left(
\Delta _{T};l_{q}\right) }+\left\Vert \varepsilon \frac{\partial ^{2}u}{%
\partial x^{2}}\right\Vert _{L_{\mathbf{p}}\left( \Delta _{T};l_{q}\right)
}+\left\Vert Au\right\Vert _{L_{\mathbf{p}}\left( \Delta _{+};E\right) }\leq
C\left\Vert f\right\Vert _{L_{\mathbf{p}}\left( \Delta _{T};E\right) }.
\end{equation*}

From the above estimate we obtain the assertion.

\begin{center}
\bigskip \textbf{9. Wentzell-Robin type mixed problem for parabolic equation
in exterior domain}
\end{center}

Consider the problem%
\begin{equation}
\frac{\partial u}{\partial t}+a\frac{\partial ^{2}u}{\partial x^{2}}+a_{1}%
\frac{\partial ^{2}u}{\partial y^{2}}+b_{1}\frac{\partial u}{\partial y}%
+cu=f\left( t,x,y\right) \text{, }  \tag{9.1}
\end{equation}%
\ \ \ 

\begin{equation}
B_{j}u=0\text{, }j=0,1,\text{ }t\in \left( 0,T\right) \text{, }x\in \sigma ,%
\text{ }y\in \left( 0,1\right) ,  \tag{9.2}
\end{equation}

\begin{equation}
u\left( 0,x,y\right) =0\text{, }x\in \sigma ,\text{ }y\in \left( 0,1\right) ,%
\text{ }  \tag{9.3}
\end{equation}

where $a=a\left( t,x,y\right) ,$ $a_{1}=a_{1}\left( t,x,y\right) ,$ $%
b_{1}=b_{1}\left( t,x,y\right) ,$ $c=c\left( t,x,y\right) $ are
complex-valued functions on $\tilde{\Omega}=\sigma \times \left( 0,1\right)
\times \left( 0,T\right) $. For $\mathbf{\tilde{p}=}\left( p,p,2\right) $
and $L_{\mathbf{\tilde{p}}}\left( \tilde{\Omega}\right) $ will denote the
space of all $\mathbf{\tilde{p}}$-summable scalar-valued\ functions with
mixed norm. Analogously, $W_{\mathbf{\tilde{p}}}^{2,1}\left( \tilde{\Omega}%
\right) $ denotes the Sobolev space with corresponding mixed norm, i.e., $W_{%
\mathbf{\tilde{p}}}^{2,1}\left( \tilde{\Omega}\right) $ denotes the space of
all functions $u\in L_{\mathbf{\tilde{p}}}\left( \tilde{\Omega}\right) $
possessing the derivatives $\frac{\partial u}{\partial t},$ $\frac{\partial
^{2}u}{\partial y^{2}},$ $\frac{\partial ^{2}u}{\partial y^{2}}\in L_{%
\mathbf{\tilde{p}}}\left( \tilde{\Omega}\right) $ with the norm 
\begin{equation*}
\ \left\Vert u\right\Vert _{W_{\mathbf{\tilde{p}}}^{2,1}\left( \tilde{\Omega}%
\right) }=\left\Vert u\right\Vert _{L_{\mathbf{\tilde{p}}}\left( \tilde{%
\Omega}\right) }+\left\Vert \frac{\partial u}{\partial t}\right\Vert _{L_{%
\mathbf{\tilde{p}}}\left( \tilde{\Omega}\right) }+\left\Vert \frac{\partial
^{2}u}{\partial x^{2}}\right\Vert _{L_{\mathbf{\tilde{p}}}\left( \tilde{%
\Omega}\right) }+\left\Vert \frac{\partial ^{2}u}{\partial y^{2}}\right\Vert
_{L_{\mathbf{\tilde{p}}}\left( \tilde{\Omega}\right) }.
\end{equation*}

\textbf{Condition 9.1 }Assume;

(1) $a\left( t,.,y\right) ,\in C\left( \bar{\sigma}\right) $, $y\in \left(
0,1\right) $ and $t\in \left( 0,T\right) ,$ $\func{Re}\omega _{k}\neq 0$ and 
$\frac{\lambda }{\omega _{k}}\in S\left( \varphi \right) $ for, $x\in \sigma
,$ $\lambda \in S\left( \varphi \right) $, $k=1,2$, $p_{k}\in \left(
1,\infty \right) ;$

(2)\ $a_{1}\left( t,x,.\right) \in W_{\infty }^{1}\left( 0,1\right) ,$ $%
a_{1}\left( t,x,.\right) \geq \delta >0,$ $b_{1}\left( t,x,.\right) ,$ $%
c\left( t,x,.\right) \in L_{\infty }\left( 0,b\right) $ for a.e. $x\in
\sigma ,$ $t\in \left( 0,T\right) ;$

\bigskip\ In this section, we present the following result:

\bigskip \bigskip \textbf{Theorem 9.1. }Suppose the Condition 9.1 hold.
Then, for $f\in L_{\mathbf{\tilde{p}}}\left( \tilde{\Omega};E\right) $
problem $\left( 9.1\right) -\left( 9.3\right) $ has a unique solution $u$\
belonging to $W_{\mathbf{\bar{p}}}^{2,1}\left( \tilde{\Omega};E\left(
A\right) ,E\right) $ and the following coercive estimate holds 
\begin{equation*}
\left\Vert \frac{\partial u}{\partial t}\right\Vert _{L_{\mathbf{\bar{p}}%
}\left( \tilde{\Omega};E\right) }+\left\Vert \frac{\partial ^{2}u}{\partial
x^{2}}\right\Vert _{L_{\mathbf{\tilde{p}}}\left( \tilde{\Omega}\right)
}+\left\Vert \frac{\partial ^{2}u}{\partial y^{2}}\right\Vert _{L_{\mathbf{%
\tilde{p}}}\left( \tilde{\Omega}\right) }+\left\Vert Au\right\Vert _{L_{%
\mathbf{\bar{p}}}\left( G_{T};E\right) }\leq C\left\Vert f\right\Vert _{L_{%
\mathbf{\bar{p}}}\left( \tilde{\Omega};E\right) }.
\end{equation*}

\ \textbf{Proof.} Let $E=L_{2}\left( 0,1\right) $. It is known $\left[ 10%
\right] $\ that $L_{2}\left( 0,1\right) $ is an $UMD$ space. Consider the
operator $A$ defined by 
\begin{equation*}
D\left( A\right) =W_{2}^{2}\left( 0,1;B_{j}u=0\right) ,\text{ }Au=a_{1}\frac{%
\partial ^{2}u}{\partial y^{2}}+b_{1}\frac{\partial u}{\partial y}+cu.
\end{equation*}

Therefore, the problem $\left( 9.1\right) -\left( 9.3\right) $ can be
rewritten in the form of $\left( 4.1\right) $, where $u\left( x\right)
=u\left( x,.\right) ,$ $f\left( x\right) =f\left( x,.\right) $\ are
functions with values in $E=L_{2}\left( 0,1\right) .$ By virtue of $\left[ 
\text{30, 31}\right] $ the operator $A$ generates analytic semigroup in $%
L_{2}\left( 0,1\right) $. Then in view of Hill-Yosida theorem (see e.g. $%
\left[ \text{28, \S\ 1.13}\right] $) this operator is sectorial in $%
L_{2}\left( 0,1\right) .$ Since all uniform bounded set in Hilbers sapace is
an $R$-bounded (see $\left[ 10\right] $ ), i.e. we get that the operator $A$
is $R$-sectorial in $L_{2}\left( 0,1\right) .$ Then from Theorem 4.1 we
obtain the assertion.

\textbf{\ }\ \textbf{References}

\begin{quote}
\ \ \ \ \ \ \ \ \ \ \ \ \ \ \ \ \ \ \ \ \ \ \ \ 
\end{quote}

\begin{enumerate}
\item Agmon S., On the eigenfunctions and on the eigenvalues of general
elliptic boundary value problems, Comm. Pure Appl. Math., 1962, 15,
119-147.\ \ \ \ \ \ \ \ \ \ \ \ \ \ \ \ \ \ \ \ \ \ \ \ \ \ \ \ \ \ \ \ \ \
\ \ \ \ \ \ \ \ \ \ \ \ \ \ \ \ \ \ \ \ \ \ \ \ \ \ \ \ \ \ \ \ \ \ \ \ 

\item Amann H., Linear and Quasi-linear Equations,1, Birkhauser, Basel, 1995.

\item Arendt W., Duelli, M., Maximal $L_{p}$- regularity for parabolic and
elliptic equations on the line, J. Evol. Equ. 2006, 6(4), 773-790.

\item Agarwal R., Bohner M., Shakhmurov V. B., Maximal regular boundary
value problems in Banach-valued weighted spaces, Boundary value problems,
(1)2005, 9-42.

\item Ashyralyev, A., Cuevas, C., and Piskarev, C., On well-posedness of
difference schemes for abstract elliptic problems in spaces, Numer. Func.
Anal. Opt., v. 29, No. 1-2, Jan. 2008, 43-65.\ 

\item Besov, O. V., Ilin, V. P., Nikolskii, S. M., Integral Representations
of Functions and Embedding Theorems, Nauka, Moscow, 1975 (in Russian).

\item Burkholder D. L., A geometrical conditions that implies the existence
certain singular integral of Banach space-valued functions, Proc. Conf.
Harmonic Analysis in Honor of Antonu Zigmund, Chicago, 1981,Wads Worth,
Belmont, 1983, 270-286.

\item Dore C. and Yakubov S., Semigroup estimates and non coercive boundary
value problems, Semigroup Forum, 2000, 60, 93-121.

\item Denk R., Hieber M., Pr\"{u}ss J., $R$-boundedness, Fourier multipliers
and problems of elliptic and parabolic type, Mem. Amer. Math. Soc. (2003),
166 (788), 1-111.

\item Favini A., Shakhmurov V., Yakubov Y., Regular boundary value problems
for complete second order elliptic differential-operator equations in UMD
Banach spaces, Semigroup Form, 2009, 79 (1), 22-54.

\item Haller R., Heck H., Noll A., Mikhlin's theorem for operator-valued
Fourier multipliers in $n$ variables, Math. Nachr. 244 (2002), 110-130.

\item Goldstain J. A., Semigroups of Linear Operators and Applications,
Oxford University Press, Oxfard, 1985.

\item Krein S. G., Linear Differential Equations in Banach
space\textquotedblright , American Mathematical Society, Providence, 1971.

\item Lunardi A., Analytic Semigroups and Optimal Regularity in Parabolic
Problems, Birkhauser, 2003.

\item Lions, J-L., Magenes, E., Nonhomogenous Boundary Value Broblems, Mir,
Moscow, 1971.

\item Sobolevskii P. E., Coerciveness inequalities for abstract parabolic
equations, Dokl. Akad. Nauk, (1964), 57(1), 27-40.

\item Shakhmurov V. B., Separable anisotropic differential operators and
applications, J. Math. Anal. Appl. 2006, 327(2), 1182-1201.

\item Shahmurov R., On strong solutions of a Robin problem modeling heat
conduction in materials with corroded boundary, Nonlinear Anal. Real World
Appl., 2011,13(1), 441-451.

\item Shahmurov R., Solution of the Dirichlet and Neumann problems for a
modified Helmholtz equation in Besov spaces on an annuals, J. Differential
Equations, 2010, 249(3), 526-550.

\item Shakhmurov V. B., Coercive boundary value problems for regular
degenerate differential-operator equations, J. Math. Anal. Appl., 2004, 292
( 2), 605-620.

\item Shakhmurov V. B., Embedding and separable differential operators in
Sobolev-Lions type spaces, Math.Notes, 2008, 84(6), 906-928.

\item Shakhmurov V. B., Separable anisotropic elliptic operators and
applications, Acta.Math. Hungar., 2011, 13(3), 208-229.

\item Shakhmurov V. B., Nonlinear abstract boundary value problems in
vector-valued function spaces and applications, Nonlinear Anal., 2006,
67(3), 745-762.

\item Triebel H., Interpolation Theory, Function Spaces, Differential
Operators, North-Holland, Amsterdam, 1978 ( in Russian).

\item Weis L, Operator-valued Fourier multiplier theorems and maximal $L_{p}$
regularity, Math. Ann., 2001, 319, 735-758.

\item Yakubov S. and Yakubov Ya., Differential-Operator Equations. Ordinary
and Partial Differential Equations, Chapman and Hall /CRC, Boca Raton, 2000.
\end{enumerate}

\end{document}